\documentclass[11pt,reqno]{amsart}
\usepackage[margin=1.15in]{geometry}
\usepackage{amsmath,amssymb,amsthm,mathrsfs}
\usepackage{tikz}
\usetikzlibrary{decorations.pathreplacing,arrows.meta}
\usepackage{xcolor}
\usepackage[colorlinks,linkcolor=blue!60!black,citecolor=blue!60!black]{hyperref}

\theoremstyle{plain}
\newtheorem{theorem}{Theorem}[section]
\newtheorem{proposition}[theorem]{Proposition}
\newtheorem{lemma}[theorem]{Lemma}
\newtheorem{corollary}[theorem]{Corollary}
\theoremstyle{definition}
\newtheorem{definition}[theorem]{Definition}
\newtheorem{remark}[theorem]{Remark}

\numberwithin{equation}{section}

\newcommand{\R}{\mathbb{R}}

\newcommand{\Abn}{\operatorname{Abn}}
\newcommand{\End}{\operatorname{End}}

\newcommand{\spn}{\operatorname{span}}
\newcommand{\D}{\Delta}
\DeclareMathOperator{\rank}{rank}

\begin{document}

\title[An equiregular counterexample to the Sard conjecture]
{An equiregular counterexample to the Sard conjecture for sub-Riemannian endpoint maps}

\author{Bangjie Wang}
\address{Bangjie Wang, \\Department of Mathematics \\Texas A\&M University\\ College Station, Texas\\ USA}
\email{bangjiewang@tamu.edu}
\urladdr{\url{https://sites.google.com/view/bangjie-wang/home?pli=1&authuser=0}}

\author{Igor Zelenko}
\address{Igor Zelenko\\
         Department of Mathematics\\
         Texas A\&M University\\
         College Station\\
         Texas \ 77843\\
         USA}
\email{zelenkotamu@tamu.edu}
\urladdr{\url{https://people.tamu.edu/~zelenkotamu/}}

\date{\today}

\begin{abstract}
Lerario, Rizzi and Tiberio disproved the sub-Riemannian Sard conjecture with a
rank-two distribution on $\R^4$, of degree of nonholonomy five at the origin, whose
endpoint map from that point has no regular values. Their distribution is neither
equiregular nor real-analytic, and the equiregular case has remained open. We construct
an equiregular rank-two distribution on an open set diffeomorphic to $\R^9$, with small
growth vector $(2,3,5,8,9)$, whose endpoint map from the origin has no regular values.
This growth vector is maximal and generic for rank-two distributions in dimension
nine, so the failure of the Sard property is not tied to any degeneracy of the growth
vector.

Singular curves of the Lerario--Rizzi--Tiberio-type structure built from the same
coefficient lift, inside the strip where our distribution is defined, to singular
curves, and the vertical displacement of a horizontal loop is an explicit strip
integral. Realising an open set of vertical displacements by singular loops reduces to a
graded moment problem, which we solve with pairs of bumps of opposite winding, an
anisotropic rescaling, and an inverse function theorem uniform over infinitely many
scales.

We also state, with proofs deferred to a planned longer version, two results about our
distribution: a description of its abnormal dynamics in the cotangent bundle, which
explains the example from a dynamical point of view, and the non-constancy of its
Tanaka symbol. The latter motivates an equiregular counterexample in dimension twelve
with constant Tanaka symbol (constant nilpotentisation), announced here and to be
proved in a separate paper. The longer version will also contain obstructions to
real-analytic and lower-dimensional versions, and a transversal variant of the
construction, which also leads to a modification of the Lerario--Rizzi--Tiberio
counterexample.

\end{abstract}

\maketitle

\section{Introduction}

Let $M$ be a smooth connected manifold, $\D\subset TM$ a smooth bracket-generating
distribution and $q_0\in M$. Writing $\Omega_{q_0}$ for the Hilbert manifold of
horizontal $W^{1,2}$ curves issuing from $q_0$ and
$\mathcal E_{q_0}:\Omega_{q_0}\to M$ for the endpoint map, the \emph{abnormal set}
is
\[
\Abn_\D(q_0):=\{\mathcal E_{q_0}(\gamma):\ \gamma\in\Omega_{q_0}\ \text{is a critical point of }\mathcal E_{q_0}\}.
\]
The Sard conjecture of Zhitomirskii and Montgomery \cite[\S10.2]{Mont} asserts that
$\Abn_\D(q_0)$ has measure zero for every smooth bracket-generating $\D$ and every
$q_0$.

In \cite{LRT} Lerario, Rizzi and Tiberio disproved this in the strongest possible
form: they produced a rank-two distribution $\D_0$ on $\R^4$ with
$\Abn_{\D_0}(0)=\R^4$, so that $\mathcal E_0$ has no regular values at all. Their
mechanism realises Kupka's counterexample \cite{Kupka} to the Morse--Sard theorem in
infinite dimensions inside a sub-Riemannian endpoint map.

Their structure is not equiregular, and it is worth being precise about how, since
the degree of nonholonomy of $\D_0$, i.e., the least $k$ with $\D_0^k=T\R^4$, depends
on the point. At the origin the small growth vector is
\[
\bigl(\dim\D_0^1,\dots,\dim\D_0^5\bigr)=(2,2,2,3,4),
\]
so the degree of nonholonomy there is five. At the points of the $x$-axis where
their coefficient is nonzero, and at the points of the second branch of the zero set
that carries the singular curves (the graph of their coefficient) where the coefficient
is nonzero, the small growth vector is $(2,2,3,4)$ and the degree of nonholonomy is
four. At all other points the small growth vector is $(2,3,4)$, the maximal possible
one, and the degree of nonholonomy is three. No other growth vector occurs. The locus of degree five is the set of points of the $x$-axis at which
the coefficient vanishes, which contains the two accumulation points of their
construction and whole intervals besides.
Equiregular structures, Carnot groups above all, are the setting of a large part of
the positive results on the Sard conjecture \cite{LMOPV,LLNPR}, and the natural question left open by \cite{LRT} is whether a
counterexample can be equiregular. We answer this question: the counterexample
constructed below is equiregular, and its growth vector is the generic one.

\subsection*{The construction}

Every structure in this paper is a \emph{graph lift} over the plane. In coordinates
$(x,y,Z)\in\R^2_{x,y}\times\R^N$ a graph lift is $\D=\spn\{X,Y\}$ with
\[
X=\partial_x+A(x,y)\cdot\partial_Z,\qquad Y=\partial_y+B(x,y)\cdot\partial_Z.
\]
A horizontal curve is determined by its base projection and initial point.
Section~\ref{sec:lifts} derives three tools from $G:=A_y-B_x$:
its derivatives determine the flag, a fixed nonzero multiplier $\mu$
with $\mu\cdot G=0$ along a base curve ensures singularity of its lift,
and integrals of $G$ compute the vertical displacements of graph loops.
The endpoint map and its derivative are also made explicit there.

The structure itself is the following. Let $a\in C_c^\infty(\R)$ satisfy
$\|a\|_\infty\le\frac12$. This coefficient is the only unknown; it is constructed in
Section~\ref{sec:assemble}. Set
\begin{equation}
\label{p_def}
p(x,y):=y\bigl(y-a(x)\bigr),\qquad w(x):=\bigl(1-x,\,x\bigr),\qquad e:=w'=(-1,1).
\end{equation}
On $\R^9$ with coordinates $(x,y,Z)$, where $Z=(z_1,z_2,v_1,\dots,v_5)\in\R^7$,
define
\begin{equation}\label{eq:frame-intro}
\begin{aligned}
X&=\partial_x
+\Bigl(\tfrac{y^3}{3}-\tfrac{a(x)y^2}{2}\Bigr)\bigl[(1-x)\partial_{z_1}+x\,\partial_{z_2}\bigr]\\
&\qquad
+(e^{y}-1)\bigl[\partial_{v_1}+e^{x}\partial_{v_2}+e^{2x}\partial_{v_3}\bigr]
+(1-e^{-y})\bigl[\partial_{v_4}+e^{x}\partial_{v_5}\bigr],\\[2pt]
Y&=\partial_y ,
\end{aligned}
\end{equation}
and $\D:=\spn\{X,Y\}$ on the open set
$M:=\{(x,y,Z):\,|y|<1\}\cong\R^9$.
We call $(x,y)$ the \emph{base} variables and $Z$ the \emph{vertical}
variables; vertical displacement means the increment of $Z$.
We use horizontal curves with $L^2$ controls that remain in $M$;
Lemma~\ref{lem:endpoint} gives the precise open control domain and
endpoint map.

For this lift, $G=A_y=(pw,h)$, where
$h=(e^y,e^{x+y},e^{2x+y},e^{-y},e^{x-y})$.
The derivatives of $h$ through order two span five dimensions and satisfy
exactly one relation, $h_{yy}=h$. Section~\ref{sec:structure} uses this
relation to prove the following uniform growth statement.

\begin{theorem}[Subsection~\ref{sec:equireg}]\label{thm:A}
$\D$ is a bracket-generating rank-two distribution on $M=\{|y|<1\}\cong\R^9$ which is
equiregular, with constant small growth vector
\begin{equation}\label{eq:growth}
\bigl(\dim\D^1,\dots,\dim\D^5\bigr)=(2,3,5,8,9).
\end{equation}
\end{theorem}

The growth vector \eqref{eq:growth} is as large as it can be. For $k\ge1$ let $d_k$
be the dimension of the free nilpotent Lie algebra of step $k$ on two generators. Then by definition of $\D^k$
\[
\dim\D^k(q)\le\min\{9,\,d_k\}\qquad\text{for every point $q$ and every }k\ge1 ,
\]
and the same bound holds, for the same reason, for every rank-two distribution on a
nine-dimensional manifold.
By
Witt's formula \cite{Reut} the homogeneous components of the free Lie algebra on two
generators have dimensions $2,1,2,3,6,9,\dots$, so that
\[
(d_1,d_2,d_3,d_4,d_5,d_6,\dots)=(2,3,5,8,14,23,\dots).
\]
For $n=9$ the right-hand sides are $2,3,5,8,9$ for $k=1,\dots,5$, which is exactly
\eqref{eq:growth}. Equality for all $k$ at a point is a nonempty open algebraic
condition on the $4$-jet of a local frame at that point (it holds, for instance, for
the distribution of Theorem~\ref{thm:A}). Hence $(2,3,5,8,9)$ is the small growth
vector of a generic rank-two distribution germ in dimension nine (compare \cite{VG}).

The factor $p$ of \eqref{p_def} carries the Lerario--Rizzi--Tiberio mechanism: its zero set is the
union of the two graphs $\{y=0\}$ and $\{y=a(x)\}$, and a horizontal curve whose base
projection runs along those graphs is singular: the multipliers
$\mu=(\nu,0)$, $\nu\in\R^2\setminus\{0\}$, annihilate $G$ there.

\begin{theorem}[Subsection~\ref{sec:singular-curves}]\label{thm:B}
Let $\nu\in\R^2\setminus\{0\}$ and let $\gamma$ be a horizontal curve of $\D$ whose
$(x,y)$-projection remains in
\[
\Sigma_\nu:=\bigl\{(x,y):\ p(x,y)\,\bigl(\nu\cdot w(x)\bigr)=0\bigr\} .
\]
Then $\gamma$ is singular, with the explicit annihilating covector
\eqref{eq:covector}.
\end{theorem}

The vertical space is seven-dimensional, so filling an open set of vertical
displacements means prescribing seven numbers. Sections~\ref{sec:gadgets}
and~\ref{sec:assemble} carry this out using the moment estimates of
Section~\ref{sec:moments}.

\begin{theorem}[Section~\ref{sec:assemble}]\label{thm:C}
There is a coefficient $a\in C_c^\infty(\R)$ with $\|a\|_\infty\le\frac12$, flat at
$0,1,2$, for which the distribution \eqref{eq:frame-intro} satisfies
\[
\Abn_\D(0)=M .
\]
Equivalently, $\mathcal E_0$ has no regular values on $M$.
\end{theorem}

We do not address here the minimizing Sard property, which concerns only the
endpoints of length-minimizing singular curves; Lerario, Rizzi and Tiberio show that
it holds for their structure \cite{LRT}.

\subsection*{The difficulty, and the shape of the answer}

A loop along the axis and back along the graph of a bump produces a
seven-dimensional vertical displacement, expressed by the moment vector
in Section~\ref{sec:moments}. These moments have leading amplitude orders
one, two and three. They cannot be prescribed independently with a single
bump: the quadratic moments are positive for every nonzero bump
\eqref{eq:Tpositive}.

We therefore pair bumps with opposite traversal directions, so their
moment contributions are subtracted. Three choices of profile and seven
independent amplitude adjustments give a \emph{gadget} of seven pairs.
After normalising the three moment blocks by their respective scales,
the moment map extends analytically to zero scales and has an invertible
derivative there. A quantitative inverse function theorem then gives
uniform local control of all seven moments
(Proposition~\ref{prop:gadget}).

Finally, we fix gadgets at infinitely many scales, with supports
accumulating at $0,1,2$ and amplitudes decaying rapidly enough that their
sum $a$ is smooth. A tail condition ensures that, in each coordinate, every
number between zero and the total increment is the sum of a subseries of the
increments.
Each selection is realised by two graph sweeps of finite energy.
Finite repetitions give arbitrary vertical displacements; appending
axis and vertical segments reaches any target in $M$.
A single multiplier certifies singularity of the whole curve, as shown
in Section~\ref{sec:assemble}.

Related nonlinear power moment problems are discussed in
\cite{SI02,SI08}; this connection is not used in the proofs.

\subsection*{Organisation}

Section~\ref{sec:lifts} develops the graph-lift calculus, and
Section~\ref{sec:structure} proves Theorems~\ref{thm:A} and~\ref{thm:B}.
Sections~\ref{sec:moments}--\ref{sec:assemble} construct the coefficient
and prove Theorem~\ref{thm:C}.
Appendices~\ref{app:symbol} and~\ref{app:cotangent} state further results
on the Tanaka symbol and abnormal dynamics without proof; their proofs
and related extensions are deferred to a longer version of the paper in preparation.

\section{Graph lifts over the plane}\label{sec:lifts}

Throughout, $Z=(Z^1,\dots,Z^N)$ are coordinates on $\R^N$ and $A,B:\R^2\to\R^N$ are
smooth. On $\R^2_{x,y}\times\R^N_Z$ consider the rank-two distribution
$\D=\spn\{X,Y\}$ with
\begin{equation}\label{eq:lift}
X=\partial_x+A(x,y)\cdot\partial_Z,\qquad
Y=\partial_y+B(x,y)\cdot\partial_Z .
\end{equation}
We call such a $\D$ a \emph{graph lift}. Everything about $\D$ is governed by the
single map
\begin{equation}\label{eq:G}
G:=A_y-B_x\ :\ \R^2\longrightarrow\R^N .
\end{equation}
The structure of \cite{LRT} is the case $N=2$, $B=0$, $A=-V$, $G=-F$.

\begin{lemma}[Flag]\label{lem:flag}
For every $k\ge2$ and every $q$,
\[
\begin{gathered}
\D^k(q)=\spn\{X,Y\}\oplus
\spn\bigl\{\partial^\alpha G(x,y)\cdot\partial_Z:\ |\alpha|\le k-2\bigr\},\\[3pt]
\dim\D^k(q)=2+\dim\spn\{\partial^\alpha G:|\alpha|\le k-2\}.
\end{gathered}
\]
\end{lemma}

\begin{proof}
Since $A,B$ do not depend on $Z$, for any smooth $K:\R^2\to\R^N$ one has
$[X,K\cdot\partial_Z]=K_x\cdot\partial_Z$ and $[Y,K\cdot\partial_Z]=K_y\cdot\partial_Z$.
Moreover $[Y,X]=(A_y-B_x)\cdot\partial_Z=G\cdot\partial_Z$. The description of
$\D^k$ follows by induction, and the sum is direct because $X,Y$ have linearly
independent $(\partial_x,\partial_y)$-components while the vertical fields have none.
\end{proof}

Fix the interval $[0,1]$ as the time domain and use the global frame: a horizontal
curve is determined by controls $(u,v)\in L^2([0,1],\R^2)$ through
\begin{equation}\label{eq:control}
\dot x=u,\qquad \dot y=v,\qquad \dot Z=uA(x,y)+vB(x,y),
\end{equation}
and $\End_{q_0}(u,v)$ denotes the resulting endpoint. The next lemma makes the
endpoint map and its derivative explicit, and identifies the open set of controls
whose trajectories stay in $M$. It is elementary, but it is what licenses the
passage from the moment computations of Section~\ref{sec:moments} to critical values
of a differentiable map, so we record it in full.

\begin{lemma}[The endpoint map]\label{lem:endpoint}
Let $A,B:\R^2\to\R^N$ be smooth and fix $q_0=(x_0,y_0,Z_0)$. For every
$c=(u,v)\in H:=L^2([0,1];\R^2)$ the system \eqref{eq:control} has a unique solution
$\gamma_c\in W^{1,2}$ with $\gamma_c(0)=q_0$, given by
\[
x_c(t)=x_0+\int_0^tu,\qquad y_c(t)=y_0+\int_0^tv,\qquad
Z_c(t)=Z_0+\int_0^t\bigl[u\,A(x_c,y_c)+v\,B(x_c,y_c)\bigr]ds .
\]
The endpoint map $\End_{q_0}:H\to\R^{N+2}$, $\End_{q_0}(c)=\gamma_c(1)$, is of class
$C^1$, and for $h=(h_1,h_2)\in H$, writing $\delta x(t)=\int_0^th_1$,
$\delta y(t)=\int_0^th_2$, $b=(x_c,y_c)$ and $\delta b=(\delta x,\delta y)$,
\begin{equation}\label{eq:variation}
D\End_{q_0}(c)\,h=\Bigl(\delta x(1),\ \delta y(1),\
\int_0^1\bigl[h_1A(b)+h_2B(b)+u\,DA(b)\delta b+v\,DB(b)\delta b\bigr]dt\Bigr).
\end{equation}
If $|y_0|<1$, the controls whose trajectories remain in $M=\{|y|<1\}$ form the open
subset
\[
\Omega_{q_0}=\Bigl\{(u,v)\in H:\ \sup_{t\in[0,1]}\Bigl|y_0+\int_0^tv\Bigr|<1\Bigr\}\subset H,
\]
and the derivative of the restriction $\End_{q_0}|_{\Omega_{q_0}}$ is again
\eqref{eq:variation}.
\end{lemma}

\begin{proof}
The primitive operator maps $L^2$ continuously into $W^{1,2}$ and into $C^0$, with
$\|\delta b\|_\infty\le\|h\|_{L^2}$. The base path is therefore continuous and
bounded, so $A,B$ are bounded together with their derivatives on its compact image,
and the displayed vertical velocity lies in $L^2$; the formulas solve
\eqref{eq:control} almost everywhere and do not depend on the representatives chosen
in the $L^2$ classes. The base equations determine $x_c,y_c$ uniquely, and the
vertical equation then determines $Z_c$ uniquely.

In a sufficiently small $L^2$-neighbourhood of $c$ all base paths lie in one compact
subset of $\R^2$. Bounding $A,B$ and their first two derivatives there, Taylor's
theorem and Cauchy--Schwarz give
\[
\bigl\|\End_{q_0}(c+h)-\End_{q_0}(c)-D\End_{q_0}(c)h\bigr\|
\le C\bigl(1+\|c\|_{L^2}\bigr)\|h\|_{L^2}^2
\]
for small $h$: the terms multiplying $u,v$ have Taylor remainders
$O(\|\delta b\|_\infty^2)$, and the terms multiplying $h_1,h_2$ leave products
bounded by $C\|h\|_{L^2}\|\delta b\|_\infty$. Hence $\End_{q_0}$ is Fréchet
differentiable with derivative \eqref{eq:variation}, and the same compact bounds make
that derivative continuous in $c$.

Finally, for an admissible control the continuity of $y_c$ on $[0,1]$ gives a
positive margin $d=1-\|y_c\|_\infty>0$; any perturbation with $\|h_2\|_{L^2}<d$
keeps the trajectory in the strip. So $\Omega_{q_0}$ is open, and restricting to an
open set does not change the derivative.
\end{proof}

Accordingly $\Abn_\D(q_0)$ is the set of $\End_{q_0}(c)$ with $c\in\Omega_{q_0}$ at
which \eqref{eq:variation} fails to be surjective. In the present construction
$N=7$ and $B\equiv0$.

\begin{lemma}[Singularity criterion, cf.\ {\cite[Lemma~2.1]{LRT}}]\label{lem:sing}
Let $\mu\in(\R^N)^*\setminus\{0\}$ and let $(u,v)$ be a control whose trajectory
$(x(t),y(t))$ satisfies
\[
\bigl(\mu\cdot G\bigr)\bigl(x(t),y(t)\bigr)\,\bigl(u(t),v(t)\bigr)=0\qquad\text{for a.e.\ }t\in[0,1];
\]
this holds in particular if $\mu\cdot G$ vanishes along the trajectory.
Then $(u,v)$ is a critical point of $\End_{q_0}$, and
\begin{equation}\label{eq:covector}
\eta:=\bigl(-\mu\cdot A(x(1),y(1)),\ -\mu\cdot B(x(1),y(1)),\ \mu\bigr)
\in T^*_{\End(u,v)}\R^{N+2}
\end{equation}
annihilates $\operatorname{Im}D_{(u,v)}\End_{q_0}$.

Conversely, if $(u,v)$ is a critical point of $\End_{q_0}$ then there is
$\mu\neq0$ with
\begin{equation}\label{eq:iff}
\bigl(\mu\cdot G\bigr)\bigl(x(t),y(t)\bigr)\,\bigl(u(t),v(t)\bigr)=0
\qquad\text{for a.e.\ }t ,
\end{equation}
that is, after discarding the intervals on which the trajectory is stationary, the
base projection lies in $\Sigma_\mu:=\{\mu\cdot G=0\}$; and every annihilating
covector is of the form \eqref{eq:covector}.

Thus criticality is equivalent to the existence of such a $\mu$. Only the first
assertion is used in the proofs of Theorems~\ref{thm:A}--\ref{thm:C}; the converse is
used in Corollary~\ref{cor:lift} and recorded for Appendix~\ref{app:cotangent}.
\end{lemma}

\begin{proof}
Write $\alpha:=\mu\cdot A$, $\beta:=\mu\cdot B$. Given a variation
$(\hat u,\hat v)\in L^2$, put $\delta x(t)=\int_0^t\hat u$, $\delta y(t)=\int_0^t\hat v$.
Differentiating \eqref{eq:control},
\[
\mu\cdot\delta Z(1)=\int_0^1(\hat u\alpha+\hat v\beta)\,dt
+\int_0^1\bigl[u(\alpha_x\delta x+\alpha_y\delta y)+v(\beta_x\delta x+\beta_y\delta y)\bigr]dt .
\]
Along the trajectory $\frac{d}{dt}\alpha=u\alpha_x+v\alpha_y$ and
$\frac{d}{dt}\beta=u\beta_x+v\beta_y$, so integrating the first integral by parts,
\[
\int_0^1\hat u\alpha\,dt=\alpha(1)\delta x(1)-\int_0^1(u\alpha_x+v\alpha_y)\delta x\,dt,
\qquad
\int_0^1\hat v\beta\,dt=\beta(1)\delta y(1)-\int_0^1(u\beta_x+v\beta_y)\delta y\,dt .
\]
Collecting the terms multiplying $\delta x$ and $\delta y$ gives
\begin{equation}\label{eq:key-identity}
\mu\cdot\delta Z(1)=\alpha(1)\delta x(1)+\beta(1)\delta y(1)
+\int_0^1(\mu\cdot G)\,\bigl(u\,\delta y-v\,\delta x\bigr)\,dt .
\end{equation}
Under the hypothesis the integrand of the last integral vanishes a.e., so
$\eta\cdot D\End(\hat u,\hat v)=-\alpha(1)\delta x(1)-\beta(1)\delta y(1)+\mu\cdot\delta Z(1)=0$.
Since $\mu\neq0$, $\eta\neq0$ and $D\End$ is not surjective.

For the converse, let $\eta=(\eta_x,\eta_y,\mu)$ annihilate the image, write
$g:=\mu\cdot G$ along the trajectory, and set $c_1:=\eta_x+\alpha(1)$,
$c_2:=\eta_y+\beta(1)$. By \eqref{eq:key-identity},
\[
\eta\cdot D\End(\hat u,\hat v)
=c_1\,\delta x(1)+c_2\,\delta y(1)+\int_0^1 g\,\bigl(u\,\delta y-v\,\delta x\bigr)\,dt ,
\]
and exchanging the order of integration in the last term this equals
\[
\int_0^1\hat u(s)\Bigl[c_1-\int_s^1 g\,v\,dt\Bigr]ds
+\int_0^1\hat v(s)\Bigl[c_2+\int_s^1 g\,u\,dt\Bigr]ds .
\]
This vanishes for every $(\hat u,\hat v)\in L^2$ if and only if both brackets vanish
for a.e.\ $s$; taking $s=1$ gives $c_1=c_2=0$, so $\eta$ has the form
\eqref{eq:covector}, and differentiating in $s$ gives $gu=gv=0$ a.e., which is
\eqref{eq:iff}. Finally $\mu\neq0$, since $\mu=0$ would force $\alpha=\beta=0$ and
hence $\eta=0$.
\end{proof}

The vertical displacement of a horizontal curve is, by \eqref{eq:control}, the line
integral $\int(A\,dx+B\,dy)$ along its base projection. We shall only ever need this
for one kind of base curve, and under the normalisation
\begin{equation}\label{eq:normalisation}
B\equiv0,\qquad A(x,y)=\int_0^yG(x,s)\,ds,
\end{equation}
for which $A(x,0)=0$: motion along $\{y=0\}$, and motion in the $y$-direction at
fixed $x$, both produce \emph{no} vertical displacement. (Any graph lift with
$B\equiv0$ can be brought to this form by a change of the vertical coordinates,
$Z\mapsto Z-\int_0^xA(s,0)\,ds$, which changes neither $\D$ nor $G$.)

\begin{lemma}[Graph-strip identity]\label{lem:flux}
Assume \eqref{eq:normalisation}. Let $b$ be smooth on $[\alpha,\beta]$ with
$b(\alpha)=b(\beta)=0$ and with graph in the strip $\{|y|<1\}$, and let $\gamma_b$ be
the closed base curve that runs along the axis from $(\alpha,0)$ to $(\beta,0)$ and
returns along the graph of $b$. Then the vertical displacement of its horizontal lift
is
\begin{equation}\label{eq:strip}
\Delta Z(\gamma_b)=-\int_\alpha^\beta A\bigl(x,b(x)\bigr)\,dx
=-\int_\alpha^\beta\!\int_0^{b(x)}G(x,y)\,dy\,dx .
\end{equation}
The inner integral is oriented, and no assumption on the sign of $b$ is needed.
Reversing $\gamma_b$ changes the sign of $\Delta Z$.
\end{lemma}

\begin{proof}
Along the lift $dZ=A\,dx$ since $B\equiv0$. On the axis segment $A(x,0)=0$, so only
the returning graph segment contributes; it is traversed from $x=\beta$ to $x=\alpha$,
whence the sign. The second form is \eqref{eq:normalisation}.
\end{proof}

This is Green's theorem for the strip between the axis and the graph, but stated
for that strip alone; the general planar formula for a closed base curve and the
signed region it bounds is true and elementary, but nothing below requires it.

\begin{remark}
Integral formulas of Stokes type have been the main tool on the \emph{positive} side
of the conjecture as well, though used in the opposite direction: Belotto da Silva
and Rifford \cite{BR} prove the Sard property on Martinet surfaces by controlling the
divergence of the vector fields generating the trace of the distribution on the
surface, so that the divergence theorem forces the abnormal set to have vanishing
volume. Here the identity is used constructively, to compute the displacements a
loop can realise rather than to bound the measure of what it reaches.
\end{remark}

Lemmas~\ref{lem:flag}--\ref{lem:flux} are the whole calculus we shall use: the flag
is read off from the derivatives of $G$, singularity from the zero sets of the
linear functionals $\mu\cdot G$, and the reachable vertical displacements from the
iterated integrals \eqref{eq:strip} of $G$ over strips bounded by those zero sets.

\section{The construction and its basic properties}\label{sec:structure}

\subsection{The structure}

Let $a\in C_c^\infty(\R)$ be real valued with $\|a\|_\infty\le\frac12$,
to be chosen in Section~\ref{sec:assemble}.
Split the vertical coordinates as $Z=(z,v)\in\R^2\times\R^5$, and set
\[
p(x,y):=y\bigl(y-a(x)\bigr),\qquad
w(x):=(1-x,x),\qquad e:=w'=(-1,1),
\]
and
\begin{equation}\label{eq:h}
h(x,y):=\bigl(e^y,e^{x+y},e^{2x+y},e^{-y},e^{x-y}\bigr)
=\bigl(e^{\alpha_i x+\beta_i y}\bigr)_{i=1}^5,
\end{equation}
where
\[
(\alpha_1,\ldots,\alpha_5)=(0,1,2,0,1),\qquad
(\beta_1,\ldots,\beta_5)=(1,1,1,-1,-1).
\]
Define
\begin{equation}\label{eq:Gdef}
G(x,y):=\bigl(p(x,y)w(x),h(x,y)\bigr)\in\R^2\times\R^5.
\end{equation}
Its first two components correspond to $z$, and its last five to $v$.
Taking $B=0$ and $A(x,0)=0$ as in \eqref{eq:normalisation}, we obtain
\[
\begin{aligned}
A(x,y)=\int_0^yG(x,t)\,dt
=\Bigl(&\Bigl(\tfrac{y^3}{3}-\tfrac{a(x)y^2}{2}\Bigr)w(x),\\
&(e^y-1)(1,e^x,e^{2x}),\ (1-e^{-y})(1,e^x)\Bigr).
\end{aligned}
\]
Thus $X=\partial_x+A\cdot\partial_Z$ and $Y=\partial_y$ give exactly
\eqref{eq:frame-intro}, and we work with $\D=\spn\{X,Y\}$ on
$M=\R\times(-1,1)\times\R^7$.
The $z$-block of $A$ has the form of the potential $V$ of \cite{LRT}, up to sign
and with their coefficient replaced by $a$; thus $\D$ projects onto the
four-dimensional structure $\D_0$ of Lerario--Rizzi--Tiberio type built from $a$
(Corollary~\ref{cor:lift}).
We now check equiregularity by differentiating $G$, then use the
factor $p$ in its $z$-block to identify singular curves.

\subsection{Equiregularity}\label{sec:equireg}

By Lemma~\ref{lem:flag}, bracket calculations reduce to counting the
independent derivatives of $G=(pw,h)$. We first count their $h$-components.
Any linear combination that cancels these components leaves a vector
in the two $z$-coordinates, so such relations reveal the additional
directions not seen by the $h$-projection.

For a multi-index $\gamma=(\gamma_1,\gamma_2)$, write
$\partial^\gamma=\partial_x^{\gamma_1}\partial_y^{\gamma_2}$. Then
\begin{equation}\label{eq:hderiv}
\partial^\gamma h=Dm_\gamma,\qquad
D:=\operatorname{diag}\bigl(e^{\alpha_i x+\beta_i y}\bigr)_{i=1}^5,
\qquad
m_\gamma:=\bigl(\alpha_i^{\gamma_1}\beta_i^{\gamma_2}\bigr)_{i=1}^5.
\end{equation}
Here zeroth powers are $1$. Since $D$ is invertible, these derivatives
have the same ranks and relations as the constant vectors $m_\gamma$.
Write $\alpha=(\alpha_i)_{i=1}^5$, $\beta=(\beta_i)_{i=1}^5$, and
$\mathbf1=(1,\ldots,1)$, with products and powers taken coordinatewise.
In particular,
\begin{equation}\label{eq:beta2}
\beta^2=\mathbf1,\qquad h_{yy}=h.
\end{equation}

\begin{lemma}\label{lem:monomials}
For the above five exponent pairs:
\begin{enumerate}
\item[(i)] $\rank\{\mathbf1,\alpha,\beta\}=3$;
\item[(ii)] $\rank\{\mathbf1,\alpha,\beta,\alpha^2,\alpha\beta\}=5$;
\item[(iii)] the only relation, up to a scalar multiple, among the six
vectors $\mathbf1,\alpha,\beta,\alpha^2,\alpha\beta,\beta^2$ is
$\beta^2-\mathbf1=0$. Equivalently, every polynomial of degree at most
two vanishing at the five exponent pairs is a multiple of $v^2-1$
in exponent-plane coordinates $(u,v)$.
\end{enumerate}
\end{lemma}

\begin{proof}
For (i), rows $1,2,4$ of the matrix with columns
$\mathbf1,\alpha,\beta$ give
$\left(\begin{smallmatrix}1&0&1\\1&1&1\\1&0&-1\end{smallmatrix}\right)$,
whose determinant is $-2$. For (ii), the matrix is
\[
\begin{pmatrix}
1&0&1&0&0\\1&1&1&1&1\\1&2&1&4&2\\1&0&-1&0&0\\1&1&-1&1&-1
\end{pmatrix},
\]
whose determinant is $-8$. Thus the six vectors in (iii) have a
one-dimensional relation space. The identity $\beta^2=\mathbf1$
provides that relation, proving (iii).
\end{proof}

\begin{proof}[Proof of Theorem~\ref{thm:A}]
Fix a point $(x,y,Z)\in M$ and write $P:=pw$, so $G=(P,h)$.
Let $E_k$ be the span of the vectors $\partial^\gamma G$ with
$|\gamma|\le k$, evaluated at this point, and set $r_k:=\dim E_k$.
Lemma~\ref{lem:flag} gives $\dim\D^{k+2}=2+r_k$.
We compute $r_0,r_1,r_2,r_3$ in order.

\emph{Orders zero and one.}
Since $h\ne0$, we have $G\ne0$ and $r_0=1$.
The $h$-components of $G,G_x,G_y$ are independent by
Lemma~\ref{lem:monomials}(i), so $r_1=3$.
Thus $\dim\D^2=3$ and $\dim\D^3=5$.

\emph{Order two: one additional $z$-direction.}
The six vectors $G,G_x,G_y,G_{xx},G_{xy},G_{yy}$ project onto a
five-dimensional space in the $h$-coordinates. By rank--nullity,
$r_2$ is five plus the dimension of the combinations with zero
$h$-component. Lemma~\ref{lem:monomials}(iii) says that every such
combination is a multiple of
\[
G_{yy}-G=(P_{yy}-P,0).
\]
Since $p_{yy}=2$ and $w$ is independent of $y$,
\begin{equation}\label{eq:R}
P_{yy}-P=(2-p)w.
\end{equation}
On $M$, $|p|=|y||y-a(x)|<3/2$, so $2-p\ne0$; also $w(x)\ne0$.
The kernel of the $h$-projection is therefore one-dimensional.
Hence $r_2=6$ and $\dim\D^4=8$.

\emph{Order three: the remaining $z$-direction.}
The $h$-projection already fills $\R^5$. It therefore suffices to
find two independent vectors in $E_3$ with zero $h$-component.
One is $G_{yy}-G$. Differentiating $h_{yy}=h$ in $x$ gives
$h_{xyy}=h_x$, so another is $G_{xyy}-G_x$.
Using $w'=e$ gives
\[
P_{xyy}=2e,\qquad P_x=p_xw+pe,\qquad
P_{xyy}-P_x=(2-p)e-p_xw.
\]
These two $z$-vectors are independent, since
\[
\det\bigl(P_{yy}-P,P_{xyy}-P_x\bigr)
=(2-p)^2\det(w,e)=(2-p)^2\ne0.
\]
Thus $r_3=5+2=7$ and $\dim\D^5=9$.

The initial rank is two because $X,Y$ have independent base components.
Together with the preceding counts, this gives $(2,3,5,8,9)$ at every
point of $M$. Thus $\D$ is equiregular and bracket-generating of step
five. Finally, $M=\R\times(-1,1)\times\R^7$ is diffeomorphic to $\R^9$.
\end{proof}

\begin{remark}\label{rem:choice}
Write $P=pw$. The choice $h_{yy}=h$ ensures that the relation among the six
$h$-derivatives of order at most two supplies a nonzero $z$-direction.
To see why its $h_{yy}$ term matters, consider a flat meeting point
$(x_0,0)$ of $y=0$ and $y=a(x)$, where $a(x_0)=a'(x_0)=0$.
For $P=pw$, direct differentiation gives, at this point,
\[
P=P_x=P_y=P_{xx}=P_{xy}=0,\qquad P_{yy}=2w(x_0).
\]
Consequently, a relation
$\sum_{|\alpha|\le2}c_\alpha\partial^\alpha h=0$
produces the $z$-vector
$\sum_{|\alpha|\le2}c_\alpha\partial^\alpha P
=2c_{(0,2)}w(x_0)$ there.
When the $h$-derivatives have rank five, this unique relation must therefore
have $c_{(0,2)}\ne0$ to obtain $\dim\D^4=8$ at such points.
Our choice gives the stronger, everywhere nonvanishing residual
$P_{yy}-P=(2-p)w$ on $M$.

For comparison, the relation $h_y=h_{xx}$ would give
\[
P_y-P_{xx}=\bigl(2y-a+a''y\bigr)w+2a'y\,e,
\]
which vanishes at every flat meeting point. If this were the unique
relation among the six $h$-derivatives, their lifts would have rank five
there, giving $\dim\D^4=7$ rather than $8$.
The coefficient constructed in Section~\ref{sec:assemble} has flat
meeting points at $0,1,2$ and at the endpoints of its bump supports;
no assumption about the order of its other zeros is needed.
\end{remark}

\subsection{Singular curves}\label{sec:singular-curves}

\begin{proof}[Proof of Theorem~\ref{thm:B}]
Fix $\nu\in\R^2\setminus\{0\}$ and set $\mu=(\nu,0)\in(\R^7)^*$.
Then $\mu\ne0$ and, by \eqref{eq:Gdef},
\[
\mu\cdot G(x,y)=p(x,y)\bigl(\nu\cdot w(x)\bigr).
\]
If the base projection of the horizontal curve remains in
$\Sigma_\nu$, this pairing vanishes along the whole curve.
Lemma~\ref{lem:sing} therefore gives singularity and the annihilating
covector \eqref{eq:covector}, with this choice of $\mu$.
\end{proof}

\begin{corollary}[Lifting singular curves]\label{cor:lift}
Let $\D_0$ be the graph lift on $\R^4$ with coordinates $(x,y,z)$, $z\in\R^2$,
spanned by $Y=\partial_y$ and the $(\partial_x,\partial_z)$-part of $X$ in
\eqref{eq:frame-intro}, so that its curvature is $G_0=pw$. Let $(u,v)\in L^2$ be a
critical point of the endpoint map of $\D_0$ whose base trajectory stays in
$\{|y|<1\}$. Then for every initial point in $M$ with the same $(x,y)$-component, the
horizontal curve of $\D$ with control $(u,v)$ is singular.
\end{corollary}

\begin{proof}
Both structures are graph lifts with the same controls, so they have the same base
trajectory. By the converse part of Lemma~\ref{lem:sing} applied to $\D_0$ there is
$\nu\ne0$ with $\bigl(\nu\cdot(pw)\bigr)(u,v)=0$ almost everywhere along it. For
$\mu=(\nu,0)$ one has $\mu\cdot G=p\,(\nu\cdot w)=\nu\cdot(pw)$, so the first part of
Lemma~\ref{lem:sing} applies to $\D$. The curve stays in $M$ because $|y|<1$.
\end{proof}

In particular, $p=0$ on $y=0$ and $y=a(x)$, so any admissible
horizontal lift along these graphs is singular. Section~\ref{sec:assemble}
chooses $a$ so that such curves can prescribe vertical displacements,
while the growth vector remains unchanged.

\section{Bump loops and moment estimates}\label{sec:moments}

\subsection{Bump loops and strip vectors}
Recall that $a\in C_c^\infty(\R)$, with
$\|a\|_\infty\le\tfrac12$, is the coefficient defining the
distribution; its specific choice will be made in
Section~\ref{sec:assemble}.

Let $I=(\alpha,\beta)$ be a bounded open interval, and let
$b\in C_c^\infty(I)$ be a bump of $a$, meaning that
$b(x)\in\{0,a(x)\}$ for every $x$.
We associate to $b$ a closed curve $\gamma_b$ in the $(x,y)$-plane as follows:
first travel along the $x$-axis from $(\alpha,0)$ to $(\beta,0)$, and then
return along the graph $y=b(x)$:
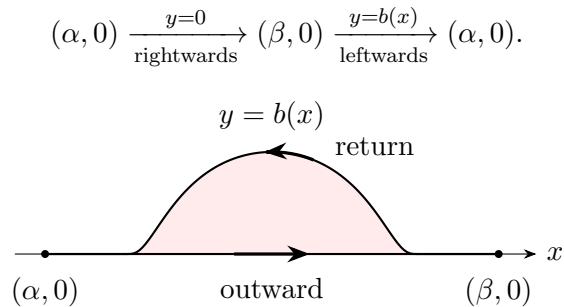
\[
(\alpha,0)
\xrightarrow[\text{rightwards}]{y=0}
(\beta,0)
\xrightarrow[\text{leftwards}]{y=b(x)}
(\alpha,0).
\]
\begin{center}
	\begin{tikzpicture}[x=1cm,y=1cm,>=Stealth]
		\fill[red!8] (1,0)
		-- plot[domain=1.01:4.99,samples=100]
		(\x,{1.35*exp(1-1/(1-((\x-3)/2)^2))})
		-- (5,0) -- cycle;
		\draw[->,thin] (-0.4,0) -- (6.5,0) node[right] {$x$};
		\draw[thick] (0,0) -- (6,0);
		\draw[thick] (0,0) -- (1,0);
		\draw[thick] plot[domain=1.01:4.99,samples=100]
		(\x,{1.35*exp(1-1/(1-((\x-3)/2)^2))});
		\draw[thick] (5,0) -- (6,0);
		\draw[very thick,->] (2.5,0) -- (3.5,0);
		\draw[very thick,->] plot[domain=3.55:2.9,samples=25]
		(\x,{1.35*exp(1-1/(1-((\x-3)/2)^2))});
		\fill (0,0) circle (1.5pt) node[below=5pt] {$(\alpha,0)$};
		\fill (6,0) circle (1.5pt) node[below=5pt] {$(\beta,0)$};
		\node[above] at (3,1.5) {$y=b(x)$};
		\node[below=5pt] at (3,0) {outward};
		\node[above right] at (3.7,1.17) {return};
	\end{tikzpicture}
	
	\small  A nonnegative bump is shown for illustration; the bumps used in the
	construction change sign.
\end{center}
Since $b$ vanishes near both endpoints of $I$, these two paths join to form
a closed, piecewise smooth curve. Since
$\|b\|_\infty\le\|a\|_\infty\le\frac12$, the curve lies in the strip $|y|<1$.
Although $\gamma_b$ is closed in the base plane, its horizontal lift need not close.
By Lemma~\ref{lem:flux}, its vertical displacement is
\begin{equation}\label{eq:cell}
\Delta Z(\gamma_b)=-\int_I A\bigl(x,b(x)\bigr)\,dx
=-\int_I\!\!\int_0^{b(x)}\!G(x,y)\,dy\,dx=:\Phi(b)\in\R^7,
\end{equation}
which we call the \emph{strip vector} of $b$; the reversed curve has displacement
$-\Phi(b)$. Note that the inner integral is oriented, so the formula remains valid
even when $b$ changes sign.

\begin{proposition}\label{prop:cell}
Let $b\in C_c^\infty(I)$ be as above, with
$b(x)\in\{0,a(x)\}$ for every $x$.
Define
\begin{equation}\label{eq:moments}
\begin{gathered}
S_m(b):=\int e^{mx}\sinh b(x)\,dx,\qquad
T_m(b):=\int e^{mx}(\cosh b(x)-1)\,dx
\qquad(m=0,1,2),\\
U(b):=\int e^{2x}\bigl(e^{b(x)}-1\bigr)\,dx
=S_2(b)+T_2(b),\\
C(b):=(C_1(b),C_2(b)):=\frac16\int b(x)^3w(x)\,dx\in\R^2.
\end{gathered}
\end{equation}
With $G$ as in \eqref{eq:Gdef}, the strip vector in the coordinate order
$(z_1,z_2,v_1,\dots,v_5)$ is
\begin{equation}\label{eq:cellv}
\begin{aligned}
\Phi(b)=\bigl(&C_1(b),C_2(b),-S_0(b)-T_0(b),-S_1(b)-T_1(b),\\
&-U(b),-S_0(b)+T_0(b),-S_1(b)+T_1(b)\bigr).
\end{aligned}
\end{equation}
\end{proposition}

\begin{proof}
For the first two coordinates, $b(x)\in\{0,a(x)\}$ implies
$a(x)b(x)^2=b(x)^3$, so
\[
\begin{aligned}
-\int_I\left(\int_0^{b(x)}p(x,y)\,dy\right)w(x)\,dx
&=-\int_I\left(\frac{b(x)^3}{3}
-\frac{a(x)b(x)^2}{2}\right)w(x)\,dx\\
&=\frac16\int_I b(x)^3w(x)\,dx=C(b).
\end{aligned}
\]
For the remaining five coordinates, since $\beta_i\in\{1,-1\}$,
\[
\int_0^{b(x)}e^{\alpha_i x+\beta_i y}\,dy
=e^{\alpha_i x}\frac{e^{\beta_i b(x)}-1}{\beta_i}
=e^{\alpha_i x}\bigl[\sinh b(x)+\beta_i(\cosh b(x)-1)\bigr].
\]
Substituting the five exponent pairs into \eqref{eq:cell} gives
\eqref{eq:cellv}.
\end{proof}

\subsection{Moment vectors}

Our goal is to prescribe the vertical displacements of the loops constructed
from bumps. The integrals in \eqref{eq:moments} define $S_m,T_m,U,C$ for every
$b\in C_c^\infty(\R)$, independently of $a$; the branch condition
$b(x)\in\{0,a(x)\}$ is needed for their identification with the strip vector
in Proposition~\ref{prop:cell}. We collect these functionals into the
\emph{moment vector}
\begin{equation}\label{eq:seven}
M(b):=\bigl(S_0(b),S_1(b),U(b);T_0(b),T_1(b);C_1(b),C_2(b)\bigr)\in\R^7.
\end{equation}
For bumps satisfying $b(x)\in\{0,a(x)\}$, Proposition~\ref{prop:cell} gives $\Phi(b)=LM(b)$,
where the fixed linear map $L:\R^7\to\R^7$ is
\[
L(s_0,s_1,u;t_0,t_1;c_1,c_2)=(c_1,c_2,-s_0-t_0,-s_1-t_1,-u,-s_0+t_0,-s_1+t_1).
\]

Since $L$ is invertible, prescribing a vertical displacement $Z_*$
is equivalent to prescribing the moment vector $L^{-1}Z_*$,
also for signed sums of bump contributions.

\subsection{Bump profiles}
For every nonzero bump $b$,
\begin{equation}\label{eq:Tpositive}
	T_m(b)=\int e^{mx}\bigl(\cosh b(x)-1\bigr)\,dx>0.
\end{equation}
Thus single bumps cannot realise arbitrary moment vectors.
We therefore use pairs with opposite traversal directions,
whose moment contributions are subtracted.
We first construct the three profiles from which these pairs will be formed.

\begin{definition}[Species]\label{def:species}
Fix once and for all three profiles $\varphi_A,\varphi_B,\varphi_C\in C_c^\infty(0,1)$
normalised by
\[
\begin{array}{llll}
\varphi_A: & \int\varphi_A=1, & &\int\varphi_A^3=0;\\[2pt]
\varphi_B: & \int\varphi_B=0, &\int\varphi_B^2=1, &\int\varphi_B^3=0;\\[2pt]
\varphi_C: & \int\varphi_C=0, & &\int\varphi_C^3=1 .
\end{array}
\]
We write $P_k^X=\int_0^1\varphi_X^k$, so $P_1^A=1$, $P_3^A=0$; $P_1^B=P_3^B=0$,
$P_2^B=1$; $P_1^C=0$, $P_3^C=1$; and $P_2^A,P_2^C>0$ automatically.
\end{definition}

\begin{lemma}\label{lem:profiles}
Such profiles exist.
\end{lemma}

\begin{proof}
Choose a nonzero, nonnegative $\eta\in C_c^\infty(0,1)$, extended by zero,
and set $m_k:=\int_0^1\eta^k>0$ for $k=1,2,3$.
For sufficiently small $h>0$, choose intervals $(s_j,s_j+h_j)$ with
pairwise disjoint closures in $(0,1)$, where $(h_1,h_2,h_3)=(h,8h,h)$.
The functions $f_j(x):=\eta((x-s_j)/h_j)$ satisfy
$\int f_j^k=h_jm_k$. Set
\begin{equation}\label{eq:profile-construction}
\varphi_A:=\frac{f_2-2f_1}{6hm_1},\qquad
\varphi_B:=\frac{f_1-f_3}{\sqrt{2hm_2}},\qquad
\varphi_C:=\frac{8f_1-f_2}{(504hm_3)^{1/3}}.
\end{equation}
All three belong to $C_c^\infty(0,1)$. Since the supports are disjoint,
\[
\int\left(\sum_{j=1}^3\lambda_jf_j\right)^k
=m_k\sum_{j=1}^3h_j\lambda_j^k,\qquad k=1,2,3.
\]
Substitution gives exactly the normalisations in Definition~\ref{def:species}.
\end{proof}

\subsection{Pairs and normalised moments}
\begin{definition}[Pairs and their moments]\label{def:pair}
	Let $X\in\{A,B,C\}$, $\ell>0$, $\kappa>0$, $x\in\R$,
	and $|\zeta|\le1$. The \emph{$X$-pair}
	$\Pi_X(x,\ell,\kappa;\zeta)$ consists of the two bumps
	\[
	b_+(s)=\kappa(1+\zeta)\,
	\varphi_X\!\left(\frac{s-x}{\ell}\right),
	\qquad
	b_-(s)=\kappa(1-\zeta)\,
	\varphi_X\!\left(\frac{s-x-\ell}{\ell}\right),
	\qquad s\in\R,
	\]
	assigned opposite traversal directions as follows.
	For $b_+$, the associated loop travels rightwards along the
	$x$-axis and returns along its graph; for $b_-$, this
	direction is reversed.
\end{definition}
\begin{remark}
	Here $\ell$ and $\kappa$ are the width and amplitude scales,
	respectively, while $\zeta$ is the adjustment parameter.
	For fixed scales and base point $x$, we vary $\zeta$ to control
	the signed moments of the pair.
	In Proposition~\ref{prop:gadget}, we combine seven pairs into a moment map
	with seven independent adjustment parameters and apply the inverse function
	theorem to realise all sufficiently small moment vectors.
	This proposition is the key local controllability result used to prove
	Theorem~\ref{thm:C}.
	To pass from Proposition~\ref{prop:gadget} to Theorem~\ref{thm:C},
	we choose a sequence of width and amplitude scales and assemble the
	corresponding gadgets into a single smooth coefficient $a$.
	Suitable selections of these gadgets then yield the required displacements.
\end{remark}
\begin{definition}[Signed moments]
The signed moments of the pair are defined by
\begin{equation}\label{eq:pair-moments}
\begin{aligned}
\Delta S_m(\Pi_X)&:=S_m(b_+)-S_m(b_-),&
\Delta T_m(\Pi_X)&:=T_m(b_+)-T_m(b_-),\\
\Delta U(\Pi_X)&:=U(b_+)-U(b_-),&
\Delta C(\Pi_X)&:=C(b_+)-C(b_-),
\end{aligned}
\end{equation}
where $m\in\{0,1,2\}$. Its moment vector is
\[
M(\Pi_X):=M(b_+)-M(b_-)=\left(
\Delta S_0(\Pi_X),\Delta S_1(\Pi_X),\Delta U(\Pi_X);
\Delta T_0(\Pi_X),\Delta T_1(\Pi_X);\Delta C(\Pi_X)\right)
.
\]
In particular, $\Delta U(\Pi_X)=\Delta S_2(\Pi_X)+\Delta T_2(\Pi_X)$.
\end{definition}

\begin{lemma}\label{lem:pair}
For $X\in\{A,B,C\}$, consider $\Pi_X(x,\ell,\kappa;\zeta)$ with
$x\in\R$, $\ell,\kappa>0$ and $|\zeta|\le1$.
The normalised moments
\[
\frac{\Delta S_m(\Pi_X)}{\ell\kappa},\qquad
\frac{\Delta U(\Pi_X)}{\ell\kappa},\qquad
\frac{\Delta T_m(\Pi_X)}{\ell\kappa^2},\qquad
\frac{\Delta C(\Pi_X)}{\ell\kappa^3},\qquad m\in\{0,1,2\},
\]
extend jointly real analytically in $(x,\zeta,\ell,\kappa)$ across
$\ell=\kappa=0$. At $(\ell,\kappa)=(0,0)$ their values are, respectively,
\[
2P_1^X e^{mx}\zeta,\qquad
2P_1^X e^{2x}\zeta,\qquad
2P_2^X e^{mx}\zeta,\qquad
P_3^X\left(\zeta+\frac13\zeta^3\right)w(x).
\]
\end{lemma}

\begin{proof}
Write $\varphi=\varphi_X$. Substituting $u=x+\ell s$ for $b_+$ and
$u=x+\ell+\ell s$ for $b_-$ in the defining integrals gives
\[
\frac{\Delta S_m(\Pi_X)}{\ell\kappa}
=e^{mx}\int_0^1 e^{m\ell s}
\left[
\frac{\sinh(\kappa(1+\zeta)\varphi(s))}{\kappa}
-e^{m\ell}\frac{\sinh(\kappa(1-\zeta)\varphi(s))}{\kappa}
\right]ds.
\]
The quotient $\sinh(\kappa a)/\kappa$ extends to an entire
function of $(\kappa,a)$, with value $a$ at $\kappa=0$.
Since $\varphi$ is bounded, the power series of
$\sinh(\kappa a)/\kappa$, after substituting
$a=(1\pm\zeta)\varphi(s)$, converges locally uniformly
in the complex parameters, uniformly for $s\in[0,1]$.
Together with the exponential factors, this justifies termwise
integration and gives a jointly analytic extension across
$\ell=\kappa=0$, whose value there is
\[
e^{mx}\int_0^1\bigl[(1+\zeta)\varphi(s)-(1-\zeta)\varphi(s)\bigr]\,ds
=e^{mx}\int_0^1 2\zeta\varphi(s)\,ds
=2P_1^Xe^{mx}\zeta.
\]

The same argument applies to $\Delta T_m/(\ell\kappa^2)$
and $\Delta U/(\ell\kappa)$ (with $m=2$ for $\Delta U$), using the entire extensions of
$(\cosh(\kappa a)-1)/\kappa^2$ and
$(e^{\kappa a}-1)/\kappa$, whose values at $\kappa=0$
are $a^2/2$ and $a$, respectively.
This gives the stated zero-scale values.

Finally, the definition $C(b)=\frac16\int b^3w$ gives directly
\[
\begin{aligned}
\frac{\Delta C(\Pi_X)}{\ell\kappa^3}
=\frac16\int_0^1\varphi(s)^3\bigl[&(1+\zeta)^3w(x+\ell s)\\
&-(1-\zeta)^3w(x+\ell+\ell s)\bigr]\,ds.
\end{aligned}
\]
Because $w$ is affine, this expression is polynomial in $(x,\zeta,\ell)$
and independent of $\kappa$. At $\ell=0$ it equals
$\frac16P_3^X[(1+\zeta)^3-(1-\zeta)^3]w(x)=P_3^X\left(\zeta+\frac13\zeta^3\right)w(x)$,
completing the proof.
\end{proof}

\section{Gadgets and local controllability}\label{sec:gadgets}

\begin{definition}[Gadget]\label{def:gadget}
Fix the macroscopic centres $0,1,2$ and $r=\frac14$. A
\emph{gadget} with data $(\ell,\kappa)$ and base points
$\mathbf x=(x_1,\dots,x_7)$ is the union of the seven pairs
\[
\Pi_A(x_1;\sigma_1),\ \Pi_A(x_2;\sigma_2),\ \Pi_A(x_3;\sigma_3),\quad
\Pi_B(x_4;\xi_1),\ \Pi_B(x_5;\xi_2),\quad
\Pi_C(x_6;\theta_1),\ \Pi_C(x_7;\theta_2),
\]
all with the same $(\ell,\kappa)$, where
\[
x_1,x_4,x_6\in(-r,r),\qquad
x_2,x_5,x_7\in(1-r,1+r),\qquad x_3\in(2-r,2+r).
\]
The base points are chosen so that all fourteen bumps have pairwise disjoint supports.
Such a choice is possible for sufficiently small $\ell$.
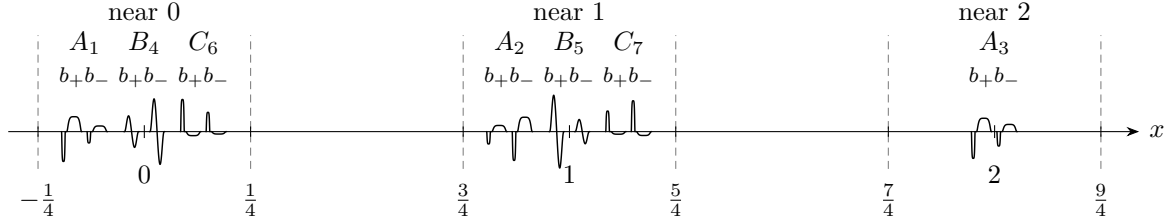
\begin{center}
\begin{tikzpicture}[x=\dimexpr\linewidth*5/14\relax,y=1.7cm,>=Stealth,font=\small]
  \draw[->] (-0.32,0) -- (2.34,0) node[right] {$x$};
  \foreach \centre in {0,1,2} {
    \draw[densely dashed,gray] (\centre-0.25,-0.29)
      -- (\centre-0.25,0.75);
    \draw[densely dashed,gray] (\centre+0.25,-0.29)
      -- (\centre+0.25,0.75);
    \draw (\centre,0.05) -- (\centre,-0.05);
    \node[below=9pt] at (\centre,0) {$\centre$};
  }
  \foreach \position/\labeltext in
    {-0.25/{-\frac14},0.25/{\frac14},
     0.75/{\frac34},1.25/{\frac54},
     1.75/{\frac74},2.25/{\frac94}} {
    \draw (\position,0.04) -- (\position,-0.04);
    \node[below=2pt] at (\position,-0.29) {$\labeltext$};
  }
  \foreach \start/\species/\idx/\imbalance/\profiletype in
    {-0.20/A/1/0.45/0,-0.06/B/4/-0.35/1,0.08/C/6/0.25/2,
      0.80/A/2/-0.40/0,0.94/B/5/0.50/1,1.08/C/7/-0.20/2,
      1.94/A/3/0.30/0} {
    \foreach \offset/\direction/\amplitudesign in {0/+/1,0.06/-/-1} {
      \ifnum\profiletype=1
        \draw[line width=0.55pt,line cap=round,line join=round]
          plot[domain=0.001:0.999,samples=160,variable=\t]
          ({\start+\offset+0.012+0.036*\t},
           {0.24*(1+\amplitudesign*\imbalance)
             *exp(1-1/(4*\t*(1-\t)))*sin(360*\t)});
      \else
        \ifnum\profiletype=0
          \def\lobes{0.02/0.18/-0.16,0.24/0.71/0.08}
        \else
          \def\lobes{0.02/0.18/0.20,0.24/0.71/-0.025}
        \fi
        \foreach \lobestart/\lobewidth/\lobeheight in \lobes {
          \draw[line width=0.55pt,line cap=round,line join=round]
            plot[domain=0.001:0.999,samples=100,variable=\t]
            ({\start+\offset+0.004+0.052*(\lobestart+\lobewidth*\t)},
             {\lobeheight*(1+\amplitudesign*\imbalance)
               *exp(-((2*\t-1)^4)/(1-(2*\t-1)^2))});
        }
      \fi
      \node[font=\scriptsize,above,inner sep=1pt]
        at (\start+\offset+0.03,0.32) {$b_{\direction}$};
    }
    \node[above,inner sep=1pt] at (\start+0.06,0.58)
      {$\species_{\idx}$};
  }
  \node[above,inner sep=1pt] at (0,0.84) {near $0$};
  \node[above,inner sep=1pt] at (1,0.84) {near $1$};
  \node[above,inner sep=1pt] at (2,0.84) {near $2$};
\end{tikzpicture}

\small
$A_1=\Pi_A(x_1;\sigma_1)$, $A_2=\Pi_A(x_2;\sigma_2)$,
$A_3=\Pi_A(x_3;\sigma_3)$;
$B_4=\Pi_B(x_4;\xi_1)$, $B_5=\Pi_B(x_5;\xi_2)$;
$C_6=\Pi_C(x_6;\theta_1)$, $C_7=\Pi_C(x_7;\theta_2)$.

Each profile combines two elementary bumps, as in
\eqref{eq:profile-construction} in the proof of Lemma~\ref{lem:profiles}.
As in \eqref{eq:profile-construction}, the $B$-profile consists of
two identical translated lobes with opposite signs, whereas
the $A$- and $C$-profiles have unequal lobe widths and heights. All use the same width parameter $\ell$; within each pair $(b_+,b_-)$, the two bumps share the same profile,
with base points separated by $\ell$ and amplitude factors
$\kappa(1+\zeta)$ and $\kappa(1-\zeta)$, respectively.
The lobe widths and vertical scales are schematic.
\end{center}

The moment vector of the gadget is the sum of moment vectors of
its pairs; we write it as
\[
\begin{aligned}
	&\qquad\Psi_{\ell,\kappa,\mathbf x}(\sigma,\xi,\theta)
	:=\sum_{j=1}^{3}
	M\bigl(\Pi_A(x_j,\ell,\kappa;\sigma_j)\bigr)\\
	&+\sum_{j=1}^{2}
	M\bigl(\Pi_B(x_{j+3},\ell,\kappa;\xi_j)\bigr)+\sum_{j=1}^{2}
	M\bigl(\Pi_C(x_{j+5},\ell,\kappa;\theta_j)\bigr)
	\in\R^7,
\end{aligned}
\]
where each pair contributes
$M(\Pi_X)=M(b_+)-M(b_-)$. The same formula defines
$\Psi_{\ell,\kappa,\mathbf x}$ for arbitrary base points.
We aim to show that, for sufficiently small $\ell$ and $\kappa$
and base points sufficiently close to the prescribed centres,
varying $(\sigma,\xi,\theta)$ near zero allows
$\Psi_{\ell,\kappa,\mathbf x}$ to attain every moment vector
in a neighbourhood of the origin in $\R^7$.

\end{definition}

\begin{lemma}[Quantitative inverse function theorem]\label{lem:qift}
Let $\|\cdot\|$ be the sup-norm on $\R^n$ and also the induced operator norm on
matrices, and let $\overline B_r$ denote the closed ball of radius $r$ about the
origin. Let $A\ge1$, let $J$ be an invertible $n\times n$ matrix with
$\|J^{-1}\|\le A$, let $\varepsilon>0$, and let $F\in C^1(\overline B_\varepsilon,\R^n)$
satisfy $\|DF(z)-J\|\le 1/(2A)$ for every $z\in\overline B_\varepsilon$.
Then $F$ is injective on $\overline B_\varepsilon$ and
$F(\overline B_\varepsilon)\supseteq
\{y:\|y-F(0)\|<\varepsilon/(2A)\}$.
\end{lemma}

\begin{proof}
This is the Banach contraction argument from the proof of the inverse function theorem
\cite[Theorem~9.24]{Rudin} (see also \cite[\S\S9.7 and~10.7]{Zorich}); it does not
depend on the choice of norm. We include it for completeness.

For $y\in\R^n$, set $\mathcal N_y(z):=z-J^{-1}(F(z)-y)$.
Its fixed points are precisely the solutions of $F(z)=y$.
Since
\[
D\mathcal N_y(z)=J^{-1}\bigl(J-DF(z)\bigr),\qquad
\|D\mathcal N_y(z)\|\le A\cdot\tfrac1{2A}=\tfrac12
\qquad(z\in\overline B_\varepsilon),
\]
the mean value inequality on the convex set $\overline B_\varepsilon$ gives
\begin{equation}\label{eq:contraction}
\|\mathcal N_y(z)-\mathcal N_y(z')\|\le\tfrac12\|z-z'\|
\qquad(z,z'\in\overline B_\varepsilon).
\end{equation}
If $\|y-F(0)\|<\varepsilon/(2A)$, then for $z\in\overline B_\varepsilon$,
\[
\|\mathcal N_y(z)\|
\le\|\mathcal N_y(z)-\mathcal N_y(0)\|+\|J^{-1}(F(0)-y)\|
<\tfrac12\|z\|+\tfrac\varepsilon2\le\varepsilon .
\]
Thus $\mathcal N_y$ is a contraction of the complete metric space
$\overline B_\varepsilon$ into itself, and the Banach fixed point theorem
provides a solution of $F(z)=y$ in $\overline B_\varepsilon$, proving the asserted
inclusion. Finally, if $F(z)=F(z')$ for some $z,z'\in\overline B_\varepsilon$,
then both are fixed points of the same map $\mathcal N_y$ with $y=F(z)$, and
\eqref{eq:contraction} forces $z=z'$, proving injectivity.
\end{proof}

\begin{proposition}[Gadget lemma]\label{prop:gadget}
There are constants $\varepsilon\in(0,\tfrac12]$, $\delta>0$ and $c_*>0$,
depending only on the fixed profiles, with the following property.
Suppose $0<\ell,\kappa<\delta$ and the base points satisfy
\(
\max\{|x_1|,|x_4|,|x_6|,|x_2-1|,|x_5-1|,|x_7-1|,|x_3-2|\}<\delta.
\)
Then the map
$\Psi_{\ell,\kappa,\mathbf x}$ restricted to $\{\|\sigma\|_\infty,\|\xi\|_\infty,\|\theta\|_\infty\le\varepsilon\}$ is injective and its image
contains the box
\(
\bigl\{y\in\R^7:\ \|y\|_\infty\le c_*\ell\kappa^3\bigr\}.
\)
\end{proposition}

\begin{proof}
Write $z=(\sigma,\xi,\theta)\in\R^7$ and
$\mathbf c=(0,1,2,0,1,0,1)$, and use sup-norms and their induced matrix norms.
Decompose $\Psi=(\Psi_{S,U};\Psi_T;\Psi_C)$ into the total signed
$(S_0,S_1,U)$, $(T_0,T_1)$ and $(C_1,C_2)$ moments, respectively.
For positive $\ell,\kappa$, define
\begin{equation}\label{eq:rescale}
\widehat\Psi_{\ell,\kappa,\mathbf x}(z)
:=\left(\frac{\Psi_{S,U}(z)}{\ell\kappa};
\frac{\Psi_T(z)}{\ell\kappa^2};
\frac{\Psi_C(z)}{\ell\kappa^3}\right).
\end{equation}
This invertible linear change of output coordinates preserves injectivity.
We will show that the normalised map is injective on
$\overline B_\varepsilon$ and its image contains $\overline B_{c_*}$.

\emph{Step 1: remove the apparent singularities in the rescaling.}
By Lemma~\ref{lem:pair}, each normalised pair moment extends jointly
analytically in its adjustment parameter, base point and scales.
Summing the seven pairs therefore extends $\widehat\Psi$ jointly
analytically near
$\{\|z\|\le1\}\times\{(\ell,\kappa,\mathbf x)=(0,0,\mathbf c)\}$.
At zero scales, this extension concerns only the moment formulas,
since the bumps themselves are not defined at $\ell=0$.

\emph{Step 2: compute the limiting map.}
Write the adjustment parameters of the seven pairs as
\[
z=(\zeta_1,\ldots,\zeta_7)
=(\sigma_1,\sigma_2,\sigma_3,\xi_1,\xi_2,\theta_1,\theta_2).
\]
Their species are $(X_1,\ldots,X_7)=(A,A,A,B,B,C,C)$,
and $c_j$ denotes the $j$th component of $\mathbf c$.
Let $\big|_{0,0,\mathbf c}$ denote evaluation of the analytic extension at
$(\ell,\kappa,\mathbf x)=(0,0,\mathbf c)$, with $z$ unchanged. Then
Lemma~\ref{lem:pair} gives
\[
\left.\frac{\Psi_{S,U}(z)}{\ell\kappa}\right|_{0,0,\mathbf c}
=
2\sum_{j=1}^7 P_1^{X_j}\zeta_j
(1,e^{c_j},e^{2c_j}),
\]
\[
\left.\frac{\Psi_T(z)}{\ell\kappa^2}\right|_{0,0,\mathbf c}
=
2\sum_{j=1}^7 P_2^{X_j}\zeta_j
(1,e^{c_j}),
\]
\[
\left.\frac{\Psi_C(z)}{\ell\kappa^3}\right|_{0,0,\mathbf c}
=
\sum_{j=1}^7 P_3^{X_j}
\left(\zeta_j+\frac13\zeta_j^3\right)w(c_j).
\]

Since $P_1^A=1$ and $P_1^B=P_1^C=0$, only the three
$A$-pairs contribute to the $(S,U)$-block.
Their parameters are $\sigma_1,\sigma_2,\sigma_3$ and their
centres are $0,1,2$, respectively. Thus
\[
\left.\frac{\Psi_{S,U}(z)}{\ell\kappa}\right|_{0,0,\mathbf c}
=
2\bigl[
\sigma_1(1,1,1)
+\sigma_2(1,e,e^2)
+\sigma_3(1,e^2,e^4)
\bigr].
\]

Similarly, all seven pairs contribute to the $T$-block.
The $A$-pairs have parameters $\sigma_1,\sigma_2,\sigma_3$ and centres
$0,1,2$; the $B$-pairs have parameters $\xi_1,\xi_2$ and centres $0,1$;
the $C$-pairs have parameters $\theta_1,\theta_2$ and centres $0,1$.
Since $P_2^B=1$, we obtain
\[
\begin{aligned}
\left.\frac{\Psi_T(z)}{\ell\kappa^2}\right|_{0,0,\mathbf c}
={}&2P_2^A\bigl[\sigma_1(1,1)+\sigma_2(1,e)+\sigma_3(1,e^2)\bigr]\\
&+2\bigl[\xi_1(1,1)+\xi_2(1,e)\bigr]\\
&+2P_2^C\bigl[\theta_1(1,1)+\theta_2(1,e)\bigr].
\end{aligned}
\]
Since $P_3^A=P_3^B=0$ and $P_3^C=1$, only the two $C$-pairs
contribute to the $C$-block. Their parameters are $\theta_1,\theta_2$
and their centres are $0,1$, respectively. Thus
\[
\begin{aligned}
\left.\frac{\Psi_C(z)}{\ell\kappa^3}\right|_{0,0,\mathbf c}
&=\left(\theta_1+\frac13\theta_1^3\right)w(0)
+\left(\theta_2+\frac13\theta_2^3\right)w(1)\\
&=\left(\theta_1+\frac13\theta_1^3,\,
\theta_2+\frac13\theta_2^3\right).
\end{aligned}
\]
In particular, $\widehat\Psi_{0,0,\mathbf c}(0)=0$.

\emph{Step 3: check invertibility at the origin.}
With columns ordered as $(\sigma_1,\sigma_2,\sigma_3,\xi_1,\xi_2,\theta_1,\theta_2)$
and rows as $(S_0,S_1,U;T_0,T_1;C_1,C_2)$, Step 2 gives
\[
J:=D_z\widehat\Psi_{0,0,\mathbf c}(0)
=\begin{pmatrix}
2&2&2&0&0&0&0\\
2&2e&2e^2&0&0&0&0\\
2&2e^2&2e^4&0&0&0&0\\
2P_2^A&2P_2^A&2P_2^A&2&2&2P_2^C&2P_2^C\\
2P_2^A&2eP_2^A&2e^2P_2^A&2&2e&2P_2^C&2eP_2^C\\
0&0&0&0&0&1&0\\
0&0&0&0&0&0&1
\end{pmatrix}.
\]
Expanding along the last two rows leaves a block lower triangular
matrix, so
\[
\det J
=2^5\det\begin{pmatrix}1&1&1\\1&e&e^2\\1&e^2&e^4\end{pmatrix}
\det\begin{pmatrix}1&1\\1&e\end{pmatrix}
=2^5(e-1)^2(e^2-1)(e^2-e)\ne0.
\]
Thus $J$ is invertible. Set $A:=\max\{1,\|J^{-1}\|\}$.

\emph{Step 4: choose uniform parameters and apply Lemma~\ref{lem:qift}.}
Set $F=\widehat\Psi_{\ell,\kappa,\mathbf x}$ and
$\varepsilon=1/(2\sqrt A)\le1/2$.
Only the last two diagonal entries of $D_z\widehat\Psi_{0,0,\mathbf c}(z)$
differ from $J$, by $\theta_1^2,\theta_2^2$. Thus
\[
\|D_z\widehat\Psi_{0,0,\mathbf c}(z)-J\|
\le\varepsilon^2=\frac1{4A}
\qquad (\|z\|\le\varepsilon).
\]
By Step 1 and compactness of $\overline B_\varepsilon$, choose
$\delta\in(0,\min\{1,r\})$ so that
\[
\sup_{\|z\|\le\varepsilon}
\|D_zF(z)-D_z\widehat\Psi_{0,0,\mathbf c}(z)\|\le\frac1{4A}
\]
whenever $0<\ell,\kappa<\delta$ and $\|\mathbf x-\mathbf c\|<\delta$.
The triangle inequality then gives $\|D_zF(z)-J\|\le1/(2A)$
on $\overline B_\varepsilon$.

To obtain an image ball centred at the origin, shrink $\delta$ further
so that $\|F(0)\|<\varepsilon/(4A)$; this follows from joint continuity
and $\widehat\Psi_{0,0,\mathbf c}(0)=0$.
Here $F(0)$ sets only the adjustment parameters to zero, not the scales.
Lemma~\ref{lem:qift} gives injectivity and an image containing the open
ball of radius $\varepsilon/(2A)$ about $F(0)$.
With $c_*:=\varepsilon/(8A)$, every $\|y\|\le c_*$ satisfies
$\|y-F(0)\|<3\varepsilon/(8A)<\varepsilon/(2A)$, hence
$F(\overline B_\varepsilon)\supseteq\overline B_{c_*}$.
All choices depend only on the fixed profiles.

\emph{Step 5: restore the original scales.}
Multiplying the three output blocks by $\ell\kappa$, $\ell\kappa^2$ and
$\ell\kappa^3$ preserves injectivity. Since $0<\kappa<\delta<1$,
each factor is at least $\ell\kappa^3$, so the resulting image contains
$\{y\in\R^7:\|y\|_\infty\le c_*\ell\kappa^3\}$.
\end{proof}

\section{Construction of the coefficient and proof of Theorem~\ref{thm:C}}\label{sec:assemble}

We first choose the scales and base points needed for the entire construction.
The following proposition makes all these choices independently of the target.

\begin{proposition}[Choice of construction parameters]\label{prop:parameters}
Fix the profiles and $r=1/4$.
Let $\varepsilon\in(0,\tfrac12]$, $\delta>0$ and $c_*>0$ be supplied by Proposition~\ref{prop:gadget},
with $\delta<\min\{1,r\}$, and set $H:=\max_{X\in\{A,B,C\}}\|\varphi_X\|_{L^\infty(0,1)}$.
There exist $K\in(2^{-1/3},1)$ and an integer $N_0\ge1$ such that, on setting
\[
\kappa_n=K^n,\qquad \ell_n=\frac{r}{50n(n+1)},\qquad
\rho_n=c_*\ell_n\kappa_n^3,\qquad
I_{c,n}=\left(c+\frac{r}{n+1},c+\frac rn\right),
\]
with $c\in\{0,1,2\}$, the following hold.
\begin{enumerate}
\item\label{item:parameters-small}
For every $n\ge N_0$,
\(
0<\ell_n,\kappa_n<\delta,\quad r/n<\delta,\quad
(1+\varepsilon)H\kappa_n\le\tfrac12.
\)
\item\label{item:parameters-place}
For each $n\ge N_0$ one can choose base points
$\mathbf x^{i,n}=(x_1^{i,n},\ldots,x_7^{i,n})$, $1\le i\le7$,
such that
\[
\begin{aligned}
\relax[x_j^{i,n},x_j^{i,n}+2\ell_n]
&\subset I_{0,n}, &&j\in\{1,4,6\},\\
\relax[x_j^{i,n},x_j^{i,n}+2\ell_n]
&\subset I_{1,n}, &&j\in\{2,5,7\},\\
\relax[x_3^{i,n},x_3^{i,n}+2\ell_n]
&\subset I_{2,n}.
\end{aligned}
\]
All these closed pair intervals are pairwise disjoint over all $i,j,n$.
In particular, all base points are within $\delta$ of their prescribed centres.
\item\label{item:parameters-decay}
For every pair of nonnegative integers $k,N$, $n^N\kappa_n\ell_n^{-k}\longrightarrow0$
as $n\to\infty$. Consequently, every bump
$b(x)=\kappa_n(1\pm\zeta)\varphi_X((x-x_0)/\ell_n)$ satisfies
\[
n^N\|D^k b\|_\infty \le C_k n^N\kappa_n\ell_n^{-k}\longrightarrow0,
\qquad C_k:=(1+\varepsilon)\max_{X\in\{A,B,C\}}\|\varphi_X^{(k)}\|_{L^\infty(0,1)},
\]
uniformly for $x_0\in\R$, $X\in\{A,B,C\}$ and $|\zeta|\le\varepsilon$.

\item\label{item:parameters-sums}
The numerical series
\[
\sum_{n\ge N_0}\ell_n\kappa_n^p\quad(p=1,2,3)
\]
are finite.
\item\label{item:parameters-tail}
The increments satisfy
\[
0<\Sigma:=\sum_{n\ge N_0}\rho_n<\infty,\qquad
\rho_n\le\sum_{k>n}\rho_k\quad\text{for every }n\ge N_0.
\]
\end{enumerate}
All these choices can be made before choosing any gadget adjustment
parameters or any target point.
\end{proposition}

\begin{proof}
Choose any $K\in(2^{-1/3},1)$, for example $K=9/10$, and define the
sequences in the statement for all $n\ge1$. We verify the five items
in order and then choose a common starting index $N_0$.

\begin{enumerate}
\item \emph{Smallness.}
Since $\ell_n,\kappa_n,r/n\to0$, all inequalities in
item~\ref{item:parameters-small} hold for $n\ge N_{\mathrm{small}}$
with a sufficiently large integer $N_{\mathrm{small}}\ge1$.

\item \emph{Placement of the pairs.}
The intervals $I_{c,n}$ are the intervals $I_n$ of \cite[\S3.1]{LRT}, with the same
$r=\frac14$, translated to the centres $c$; here each of them carries many narrow
pairs instead of a single bump.
At scale $n$, the seven gadgets require $21$ pairs near each of $0$ and $1$,
and $7$ near $2$. Each interval $I_{c,n}$ has length $|I_{c,n}|=\frac{r}{n(n+1)}=50\ell_n$,
whereas the pair intervals occupy at most $21(2\ell_n)=42\ell_n$.
The remaining length allows positive gaps between the closed pair intervals
and at both endpoints of $I_{c,n}$. Choose their left endpoints as the
base points, assigning pairs $1,4,6$ of each gadget near $0$,
pairs $2,5,7$ near $1$, and pair $3$ near $2$.
The open intervals $I_{c,n}$ are disjoint for distinct scales at a fixed
centre, and lie in disjoint neighbourhoods for distinct centres.
Hence all chosen closed pair intervals are pairwise disjoint.
For $n\ge N_{\mathrm{small}}$, each base point is at distance less than
$r/n<\delta$ from its prescribed centre.

\item \emph{Decay estimates.}
For any fixed nonnegative integers $k,N$,
\[
n^N\kappa_n\ell_n^{-k}
=\left(\frac{50}{r}\right)^k K^n n^N[n(n+1)]^k
\longrightarrow0,
\]
because exponential decay dominates every polynomial.
No further restriction on $N_0$ depending on $k$ or $N$ is needed.

\item \emph{Summability.}
Since $\ell_n\kappa_n^p\le C K^{pn}$ for $p=1,2,3$ and $0<K<1$,
all three series converge.

\item \emph{The tail condition.}
The numbers $\rho_n$ are positive and summable by item~\ref{item:parameters-sums}
with $p=3$. Moreover,
\[
\frac{\sum_{k>n}\rho_k}{\rho_n}
=\sum_{j\ge1}K^{3j}\frac{n(n+1)}{(n+j)(n+j+1)}
\longrightarrow\sum_{j\ge1}K^{3j}
=\frac{K^3}{1-K^3}>1.
\]
Indeed, for each fixed $j$ the fraction tends to $1$ and is bounded by $1$,
so dominated convergence applies with summable majorant $K^{3j}$.
Consequently, there is an integer $N_{\mathrm{tail}}\ge1$ such that
$\rho_n\le\sum_{k>n}\rho_k$ for every $n\ge N_{\mathrm{tail}}$.
\end{enumerate}

Finally, take $N_0=\max\{N_{\mathrm{small}},N_{\mathrm{tail}}\}$ and retain
the base points chosen above for $n\ge N_0$.
All five properties hold, independently of the adjustment parameters
and target point.
\end{proof}

Fix the parameters and base points supplied by Proposition~\ref{prop:parameters}.
For each $i\in\{1,\ldots,7\}$ and $n\ge N_0$, Proposition~\ref{prop:gadget}
applies at $\mathbf x^{i,n}$ and directly provides a gadget $\mathcal G_{i,n}$,
with adjustment parameters in $\overline B_\varepsilon$, such that
\begin{equation}\label{eq:gadgetvalue}
M(\mathcal G_{i,n})=\rho_n\mathbf e_i.
\end{equation}
There are $98$ bumps at each scale, all with disjoint supports. Define
\[
a_{i,n}:=\sum_{\text{bumps of }\mathcal G_{i,n}}b,
\qquad a:=\sum_{n\ge N_0}\sum_{i=1}^7a_{i,n}.
\]
These are ordinary sums of functions; the winding signs occur only in the moment
vectors $M(\mathcal G_{i,n})=\sum_{\text{pairs}}[M(b_+)-M(b_-)]$.
The coefficient $a$ and all gadgets are now fixed independently of the target.

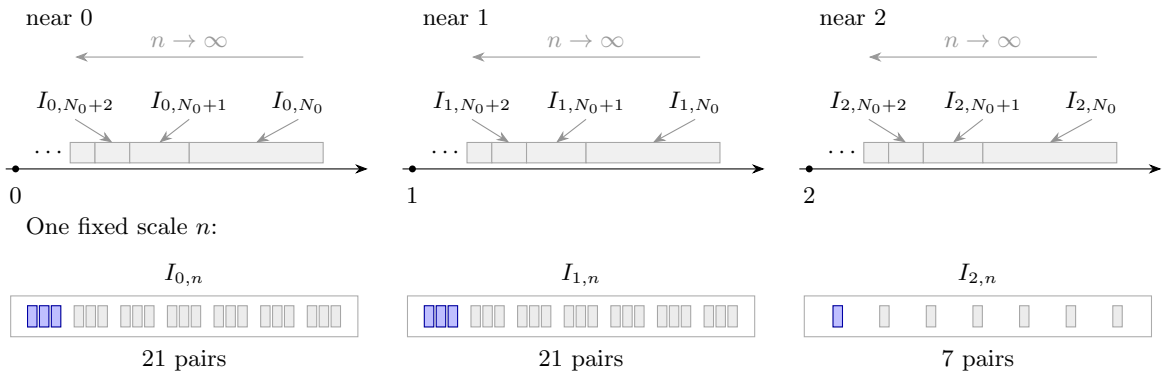
\begin{center}
\begin{tikzpicture}[x={\dimexpr\linewidth/12\relax},y=1cm,>=Stealth,font=\footnotesize]
\path[use as bounding box] (-0.2,-2.65) rectangle (11.8,2.25);
\foreach \panel/\centre in {0/0,4/1,8/2} {
  \begin{scope}[shift={(\panel,0)}]
    \node[anchor=west] at (0,2.0) {near $\centre$};
    \draw[->] (-0.1,0) -- (3.55,0);
    \fill (0,0) circle (1.3pt) node[below=3pt] {$\centre$};
    \node at (0.35,0.22) {$\cdots$};
    \foreach \leftend/\rightend in {0.55/0.8,0.8/1.15,1.15/1.75,1.75/3.1} {
      \draw[gray!70,fill=gray!12] (\leftend,0.08) rectangle (\rightend,0.35);
    }
    \node[inner sep=2pt] (thirdinterval) at (0.6,0.88) {$I_{\centre,N_0+2}$};
    \node[inner sep=2pt] (nextinterval) at (1.75,0.88) {$I_{\centre,N_0+1}$};
    \node[inner sep=2pt] (firstinterval) at (2.85,0.88) {$I_{\centre,N_0}$};
    \draw[->,gray!85] (thirdinterval.south) -- (0.975,0.35);
    \draw[->,gray!85] (nextinterval.south) -- (1.45,0.35);
    \draw[->,gray!85] (firstinterval.south) -- (2.425,0.35);
    \draw[->,gray!80] (2.9,1.48) -- (0.6,1.48)
      node[midway,above] {$n\to\infty$};
  \end{scope}
}
\node[anchor=west] at (0,-0.75) {One fixed scale $n$:};
\begin{scope}[yshift=-0.5cm]
\foreach \panel/\centre/\paircount in {0/0/21,4/1/21,8/2/7} {
  \begin{scope}[shift={(\panel,0)}]
    \draw[gray!65] (-0.05,-1.72) rectangle (3.45,-1.18);
    \node[above] at (1.7,-1.18) {$I_{\centre,n}$};
    \foreach \group in {0,...,6} {
      \ifnum\centre<2
        \foreach \slot in {0,1,2} {
          \draw[gray!70,fill=gray!15]
            ({0.12+0.47*\group+0.12*\slot},-1.59)
            rectangle ++(0.095,0.28);
        }
      \else
        \draw[gray!70,fill=gray!15]
          ({0.24+0.47*\group},-1.59) rectangle ++(0.095,0.28);
      \fi
    }
    \ifnum\centre<2
      \foreach \slot in {0,1,2} {
        \draw[blue!65!black,fill=blue!25]
          ({0.12+0.12*\slot},-1.59) rectangle ++(0.095,0.28);
      }
    \else
      \draw[blue!65!black,fill=blue!25]
        (0.24,-1.59) rectangle ++(0.095,0.28);
    \fi
    \node[below] at (1.7,-1.76) {$\paircount$ pairs};
  \end{scope}
}
\end{scope}
\end{tikzpicture}
\par\smallskip
{\footnotesize\color{blue!65!black}
The blue pairs belong to
one gadget $\mathcal G_{1,n}$ and together define $a_{1,n}$. The remaining grey pairs are grouped into six further gadgets,
whose bump sums define $a_{2,n},\ldots,a_{7,n}$, respectively.\par}
\smallskip
\small
Schematic arrangement (not to scale).
The intervals $I_{c,n}$ accumulate at $c=0,1,2$ from the right;
each small rectangle in the lower panel represents a pair of bumps. Each scale contains seven gadgets, hence $98$ bumps.
Smoothness of their sum is proved in Lemma~\ref{lem:smooth}.
\end{center}

\begin{lemma}\label{lem:smooth}
Let $\mathscr B$ denote the family of all individual bumps
in the fixed gadgets $\mathcal G_{i,n}$,
$i=1,\ldots,7$, $n\ge N_0$.
For every subfamily $E\subseteq\mathscr B$, finite or infinite,
the sum
\(
b_E:=\sum_{b\in E}b
\)
belongs to $C_c^\infty(\R)$, is flat at $0,1,2$,
and satisfies $\|b_E\|_\infty\le\tfrac12$.
In particular, taking $E=\mathscr B$ gives $b_E=a$.
Uniformly over all selections and for every integer $k\ge0$,
\begin{equation}\label{eq:subfamily}
\sup_{x\in\bigcup_{c\in\{0,1,2\}} I_{c,n}}|D^kb_E(x)|
\le C_k\kappa_n\ell_n^{-k}.
\end{equation}
The moments satisfy
\begin{equation}\label{eq:momabs}
\begin{aligned}
\sum_{b\in E}\bigl(|S_0(b)|+|S_1(b)|+|U(b)|\bigr)
&\le C\sum_{n\ge N_0}\ell_n\kappa_n<\infty,\\
\sum_{b\in E}\bigl(|T_0(b)|+|T_1(b)|\bigr)
&\le C\sum_{n\ge N_0}\ell_n\kappa_n^2<\infty,\\
\sum_{b\in E}\|C(b)\|_\infty
&\le C\sum_{n\ge N_0}\ell_n\kappa_n^3<\infty.
\end{aligned}
\end{equation}
Moreover,
\[
M(b_E)=\sum_{b\in E}M(b),\qquad
\Phi(b_E)=\sum_{b\in E}\Phi(b),
\]
with both series absolutely convergent.
\end{lemma}

\begin{proof}
Scaling gives $\|D^kb\|_\infty\le C_k\kappa_n\ell_n^{-k}$ for each
scale-$n$ bump. Disjointness gives \eqref{eq:subfamily} for $b_E$.
The family is locally finite away from $0,1,2$, so $b_E$ is smooth there.
On $I_{c,n}$, $|x-c|\asymp n^{-1}$; thus
Proposition~\ref{prop:parameters}(\ref{item:parameters-decay}) implies $|D^kb_E(x)|=o(|x-c|^N)\quad(x\to c)$
for every $k,N\ge0$. Outside the bump supports these derivatives vanish.
Extend each $D^k b_E$ by zero at the centres. The case $N=0$
gives continuity, while $N=1$ gives
$D^k b_E(x)/(x-c)\to0$, so the extended function has derivative zero at $c$.
These extensions therefore remain successive derivatives of $b_E$,
proving smoothness and flatness at $0,1,2$.
All supports lie in $[0,2+r]$, and the amplitude bound in
Proposition~\ref{prop:parameters}(\ref{item:parameters-small}) gives
$\|b_E\|_\infty\le1/2$.

The weights are bounded on $[0,2+r]$. For $|u|\le1/2$,
$|\sinh u|+|e^u-1|\le C|u|$ and $|\cosh u-1|\le Cu^2$.
Consequently, for each scale-$n$ bump and $m\in\{0,1,2\}$,
\[
\begin{aligned}
\int e^{mx}|\sinh b|\,dx+\int e^{2x}|e^b-1|\,dx
&\le C\ell_n\kappa_n,\\
\int e^{mx}|\cosh b-1|\,dx&\le C\ell_n\kappa_n^2,\\
\frac16\int |b|^3\|w(x)\|_\infty\,dx&\le C\ell_n\kappa_n^3.
\end{aligned}
\]
Summing these absolute integral bounds and applying
Proposition~\ref{prop:parameters}(\ref{item:parameters-sums}) proves
\eqref{eq:momabs}. Disjointness and $f(0)=0$ give
$f(b_E)=\sum_{b\in E}f(b)$ for
$f(u)=\sinh u,\cosh u-1,e^u-1,u^3$.
The same bounds justify termwise integration, proving additivity of
$M$. Finally, $b_E(x)\in\{0,a(x)\}$, so
Proposition~\ref{prop:cell} gives $\Phi(b_E)=LM(b_E)$
and hence the asserted absolute convergence and additivity of $\Phi$.
\end{proof}

\begin{lemma}[Subsums]\label{lem:subsum}
Let $\rho_n>0$ $(n\ge N_0)$ with $\Sigma=\sum_n\rho_n<\infty$ and
$\rho_n\le\sum_{k>n}\rho_k$ for every $n$. Then every $s\in[0,\Sigma]$ has the form
$s=\sum_nd_n\rho_n$ with $d_n\in\{0,1\}$.
\end{lemma}

\begin{proof}
Choose the digits greedily: put $s_{N_0}:=s$, and for $n\ge N_0$ let $d_n:=1$ if
$\rho_n\le s_n$ and $d_n:=0$ otherwise, and $s_{n+1}:=s_n-d_n\rho_n$. We claim
$0\le s_n\le R_n:=\sum_{k\ge n}\rho_k$ for all $n$. This holds for $n=N_0$. If it
holds for $n$ and $d_n=1$, then $0\le s_n-\rho_n\le R_n-\rho_n=R_{n+1}$; if $d_n=0$,
then $s_{n+1}=s_n<\rho_n\le R_{n+1}$ by the hypothesis. Since $R_n\to0$, we get
$s-\sum_{k<n}d_k\rho_k=s_n\to0$.
\end{proof}

\begin{proof}[Proof of Theorem~\ref{thm:C}]
Use the fixed coefficient $a$ constructed above. Recall that
$\Sigma=\sum_{n\ge N_0}\rho_n$ was defined in
Proposition~\ref{prop:parameters}(\ref{item:parameters-tail}).
We first realise moment vectors in $Q:=[0,\Sigma]^7$ by selecting gadgets.
By the tail condition in that proposition and Lemma~\ref{lem:subsum},
for every $s\in[0,\Sigma]$ there exist $d_n\in\{0,1\}$ such that
$s=\sum_{n\ge N_0}d_n\rho_n$.
Applying this coordinatewise, every $q\in Q$ admits a representation
\[
q=\sum_{n\ge N_0}\sum_{i=1}^7d_{i,n}\rho_n\mathbf e_i,
\qquad d_{i,n}\in\{0,1\}.
\]

Select precisely those gadgets $\mathcal G_{i,n}$ for which
$d_{i,n}=1$. From each selected gadget, include all seven
winding-$+1$ bumps in $f_+$ and all seven winding-$-1$ bumps
in $f_-$. Thus
\[
f_\pm:=\sum_{n\ge N_0}\sum_{i=1}^7 d_{i,n}
\sum_{\text{pairs in }\mathcal G_{i,n}} b_\pm.
\]
Lemma~\ref{lem:smooth} gives $f_\pm\in C_c^\infty(\R)$ with $\|f_\pm\|_\infty\le1/2$ and
$f_\pm(x)\in\{0,a(x)\}$. Both vanish at $x=0$ and $x=3$.
Consider the following two closed base curves, concatenated at $(0,0)$:
\begin{enumerate}
\item travel along the axis from $(0,0)$ to $(3,0)$ and return along
$y=f_+(x)$;
\item travel along $y=f_-(x)$ from $(0,0)$ to $(3,0)$ and return along
 the axis.
\end{enumerate}
This selection of bumps, by following either their graphs or the axis, is the binary
mechanism of \cite[Proposition~4.1]{LRT}, with the binary digits replaced by the
digits $d_{i,n}$ provided by Lemma~\ref{lem:subsum}.
Their horizontal lift has total vertical increment
\[
\begin{aligned}
\Delta Z
&=\Phi(f_+)-\Phi(f_-)\\
&=L\left(\sum_{n\ge N_0}\sum_{i=1}^7d_{i,n}
\sum_{\text{pairs in }\mathcal G_{i,n}}
[M(b_+)-M(b_-)]\right)\\
&=L\left(\sum_{n\ge N_0}\sum_{i=1}^7d_{i,n}\rho_n\mathbf e_i\right)
=Lq.
\end{aligned}
\]
Here $L$ is the fixed invertible map from moment to displacement
coordinates. We use the additivity and absolute convergence from
Lemma~\ref{lem:smooth}, together with \eqref{eq:gadgetvalue}
and the chosen representation of $q$.
The base lies in $\{p=0\}$. Parametrising each graph with unit horizontal
speed gives controls $u=\dot x,\quad v=f_\pm'(x)\dot x,\quad |\dot x|=1$,
so the two loops have finite $L^2$ energy since $f_\pm$ are smooth on $[0,3]$.

We next obtain any vertical increment $Z_*\in\R^7$ using finitely many
such loops. Choose an integer $N\ge1$ with
$\|L^{-1}Z_*\|_\infty\le N\Sigma$, and write
$L^{-1}Z_*/N=q_+-q_-$ with $q_+,q_-\in Q$.
Concatenate the loop for $q_+$ with the reverse of the loop for $q_-$,
and repeat this concatenation $N$ times.
Since the vector fields are independent of $Z$, increments add,
so the resulting lift ends at $(0,0,Z_*)$.
Its base still lies in $\{p=0\}$; finite concatenation preserves
finite $L^2$ energy.

Finally, fix an arbitrary target $(x_*,y_*,Z_*)\in M$ and append the
broken line
\[
(0,0)\longrightarrow(x_*,0)\longrightarrow(x_*,y_*).
\]
This is the broken line of \cite[Lemma~4.3]{LRT}, and the multiplier $\nu_*$ used
below is the one used there.
It contributes no vertical displacement: $A(x,0)=0$ on the first segment,
and $dx=0$ and $B=0$ on the second. The lift therefore reaches the target.
To verify singularity of the entire curve, take
$\nu_*=(x_*,x_*-1)\ne0$. Then
\begin{equation}\label{eq:nuw}
\nu_*\cdot w(x)=x_*-x.
\end{equation}
Thus every preceding loop and the appended broken line lie in the
single set $\Sigma_{\nu_*}=\{p(x,y)(\nu_*\cdot w(x))=0\}$.
Theorem~\ref{thm:B} applies to the whole horizontal curve.
It stays in $M$, because the loops have $|y|\le1/2$ and the final
segment has $|y|\le|y_*|<1$.
An affine rescaling of its finite time interval to $[0,1]$ preserves
finite $L^2$ energy and the singularity criterion. By
Lemma~\ref{lem:sing}, its control is a critical point of $\End_0$.
Since the target was arbitrary, $\Abn_\D(0)=M$.
\end{proof}

\appendix

\section{The Tanaka symbol of \texorpdfstring{$\D$}{D} is not constant}\label{app:symbol}

This appendix states, without proof, what the nilpotent approximation of $\D$ looks
like. Proofs are in the longer version of the paper in preparation. Nothing in
Sections~\ref{sec:lifts}--\ref{sec:assemble} depends on anything here.

For an equiregular distribution the \emph{Tanaka symbol} at $q$, introduced in
\cite{tan1}, is the graded nilpotent Lie algebra
\[
\mathfrak m(q)=\bigoplus_{i\ge1}\mathfrak g_{-i}(q),\qquad
\mathfrak g_{-1}(q)=\D(q),\qquad
\mathfrak g_{-i}(q)=\D^i(q)/\D^{i-1}(q)\ (i\ge2),
\]
with the bracket induced by the Lie bracket of vector fields; it is the
nilpotentisation of $\D$ at $q$. For our $\D$ the layer dimensions are
$(2,1,2,3,1)$ at every point, by Theorem~\ref{thm:A}.

\begin{proposition}[Graph lifts are metabelian]\label{app:metab}
Let $\D$ be any graph lift, equiregular of step five near $q$. Then
$V:=\mathfrak g_{-2}\oplus\cdots\oplus\mathfrak g_{-5}$ is an abelian ideal of
$\mathfrak m(q)$, on which $X$ and $Y$ act as $\partial_x$ and $\partial_y$. Hence
$[\mathfrak m,\mathfrak m]=V$ is abelian, so its Tanaka symbol at $q$ is metabelian.
\end{proposition}

For growth $(2,3,5,8,9)$ a metabelian symbol is determined, up to isomorphism, by a
binary cubic form on $\mathfrak g_{-1}$ taken up to a nonzero factor (namely the
class of $\operatorname{ad}_\ell^3G$ in $\mathfrak g_{-5}$), and two such symbols
are isomorphic exactly when the configurations of the roots of their cubics in
$\mathbb P(\mathfrak g_{-1})\cong\mathbb{RP}^1$ agree. Over $\R$ a nonzero binary
cubic has one of four root configurations: three distinct real roots, one real and
two complex conjugate roots, a double root and a simple root, or a triple root. (The
identically vanishing cubic is a fifth, degenerate case; it and the non-metabelian
symbols, which form a further two-parameter family, do not occur for graph lifts.)

\begin{proposition}[The cubic of $\D$]\label{app:cubic}
Write $\ell=sX+tY$. At $q=(x,y,Z)\in M$ the cubic is, up to a nonzero factor,
\[
F(s,t)=s\bigl(E_{30}s^2+3E_{21}st+3E_{12}t^2\bigr),
\]
\[
E_{12}=2-p,\qquad E_{21}=-2a'(1-y)+p+a-2y,\qquad E_{30}=2p+6a'y-3a''y .
\]
In particular $\ell=Y$ is a root at every point, and the type is decided by the sign
of the discriminant of the quadratic factor, $\mathcal D:=9E_{21}^2-12E_{30}E_{12}$.
Since $E_{12}=2-p\ge\frac12$ on $M$, the cubic is never zero and never a perfect
cube, so only three of the four configurations can occur.
\end{proposition}

\begin{theorem}[Non-constancy]\label{app:jumps}
The Tanaka symbol of $\D$ is not constant. On the axis $\{y=0\}$ one has
$\mathcal D=9(a-2a')^2\ge0$, so the axis carries the double-root type exactly on the
critical set $\{a=2a'\}=\{(e^{-x/2}a)'=0\}$ and the three-real-root type elsewhere
on the axis. Off the axis, at points where $a=a'=a''=0$,
$\mathcal D=y^2(33y^2-36y-12)$, which is negative for $0<|y|$ small. In particular
the symbol is of double-root type at the origin and of one-real-two-complex type at
$(0,y)$ for all sufficiently small $y\ne0$, so it is not constant on any
neighbourhood of the origin.
\end{theorem}

\begin{remark}[Two cautions]\label{app:cautions}
Neither statement should be read as a description of a whole stratum. The
double-root locus on the axis is not exhausted by the origin, the accumulation points
and the flat stretch: on each nonzero bump the function $e^{-x/2}a$ vanishes at both
endpoints and not identically, so it has an interior critical point, and there
$a=2a'$ with $a\ne0$. And $\mathcal D$ is not of one sign off the axis: on
$\{a=a'=a''=0\}$ it vanishes at $y_-=(6-4\sqrt5)/11=-0.2676\ldots$ and is positive
for $y\in(-1,y_-)$. We have not determined the partition of $M$ into types, and
nothing above needs it.
\end{remark}

\begin{corollary}[No flat model explains it]\label{app:flat}
By Proposition~\ref{app:metab} the Tanaka symbol $\mathfrak m(q)$ of $\D$ is
metabelian at every point $q$. Hence the flat structure with symbol $\mathfrak m(q)$,
that is, the left-invariant distribution generated by $\mathfrak g_{-1}(q)$ on the
simply connected Carnot group with Lie algebra $\mathfrak m(q)$, is a rank-two
polarization of a metabelian Lie group in the sense of \cite{LLNPR}, and such
polarizations satisfy the Sard property by \cite{LLNPR}. Therefore no nilpotent
approximation of $\D$ is itself a counterexample: the failure of the Sard property in
Theorem~\ref{thm:C} is genuinely non-homogeneous, and is not inherited from a Carnot
group.
\end{corollary}

The reason the symbol varies is structural. The coefficient $a$ enters through
$p=y(y-a)$, which contributes at weight zero, and the weight-zero part of the
structure is exactly what the symbol records. To obtain a constant symbol one must
push $a$ into structure functions of positive weight, at the cost of dimension. In a
paper in preparation \cite{ZW12} we carry this out with an explicit graph lift on
$\R^{12}$ with constant growth vector $(2,3,5,8,12)$, whose symbol is then, by
Proposition~\ref{app:metab}, necessarily the free metabelian rank-two step-five
algebra, and we solve the corresponding ten-dimensional moment problem. We announce the result: there is an
equiregular counterexample to the Sard conjecture in dimension $12$ with constant
Tanaka symbol. Its proof is not given here.

\section{The abnormal dynamics in \texorpdfstring{$T^*M$}{T*M}}\label{app:cotangent}

This appendix describes, again without proof, the construction from the point of view
of dynamical systems, through its abnormal extremals in $T^*M$, in the spirit of the
analysis of abnormal extremals of rank-two distributions in \cite{ZZ95,Z99}; a closely
related planar dynamics is the main tool of \cite{BCJPS}, where it is used to prove
$C^1$ regularity of abnormal minimisers. In our view this is the most transparent
explanation of the mechanism, although the analytic argument of
Sections~\ref{sec:moments}--\ref{sec:assemble} cannot be avoided in the proof. Nothing
above depends on it.

\subsection*{Reduction to planar Hamiltonian systems}
Use canonical coordinates $(x,y,Z;\xi,\eta,\mu)$ on $T^*M$. A covector annihilates $\D$ if
and only if $\eta=0$ and $\xi=-\mu\cdot A(x,y)$, and it annihilates $\D^2$ if and only if,
in addition, $g_\mu(x,y):=\mu\cdot G(x,y)=0$. Since $X$ and $Y$ do not depend on $Z$, the
vertical part $\mu$ is a first integral of the abnormal flow; this is why the multiplier in
Lemma~\ref{lem:sing} is constant. Along an abnormal extremal with multiplier $\mu$ the base
velocity $(u,v)$ satisfies $u\,\partial_xg_\mu+v\,\partial_yg_\mu=0$. Hence, for each fixed
$[\mu]\in\mathbb{RP}^6$, the abnormal dynamics descends to the $(x,y)$-plane: the base
projections of the abnormal extremals with multiplier $\mu$ are the curves in the zero level
$\{g_\mu=0\}$, and where $dg_\mu\ne0$ they are, up to reparametrisation, orbits of the
Hamiltonian vector field
\[
\mathcal X_\mu:=\partial_yg_\mu\,\partial_x-\partial_xg_\mu\,\partial_y .
\]
So the characteristic line field of $(\D^2)^\perp$, which lives on a manifold of dimension
$15$, becomes a family of planar Hamiltonian systems, parametrised by $[\mu]$ and restricted
to their zero levels.

\subsection*{Singular points}
A covector in $(\D^2)^\perp$ with vertical part $\mu$ annihilates $\D^3$ if and only if
$dg_\mu=0$ at its base point. There the kernel of the symplectic form restricted to
$(\D^2)^\perp$ jumps from dimension one to dimension three, so the direction of an abnormal
extremal is no longer determined and the extremal may branch. Thus the singular points of
the abnormal dynamics are the equilibria of $\mathcal X_\mu$ on its zero level. Their
classification into hyperbolic, elliptic and nilpotent points \cite{ZZ95,Z99} becomes the
classification of the critical points of $g_\mu$. Indeed, write these equations in the planar
form $\dot h=Ah/|h|$ used in \cite{BCJPS}, with
$h:=\bigl(\langle\lambda,[X,[X,Y]]\rangle,\langle\lambda,[Y,[X,Y]]\rangle\bigr)$ and $A$
traceless. Here $h=-dg_\mu$ and $A=\pm\operatorname{Hess}(g_\mu)\,J$, where $J$ is
the rotation by a right angle. So $\det A=\det\operatorname{Hess}g_\mu$, and a singular point
is hyperbolic, elliptic or nilpotent according as $g_\mu$ has there a saddle (two branches of
the zero level cross), an extremum (an isolated point of the zero level) or a degenerate
critical point.

\subsection*{The phase portrait for the multipliers \texorpdfstring{$(\nu,0)$}{(nu,0)}}
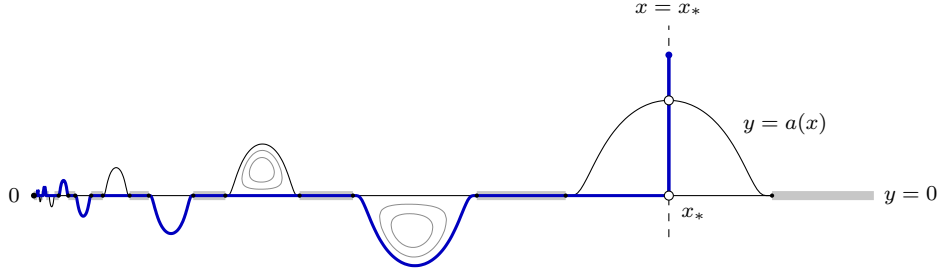
\begin{figure}[ht]
\centering
\begin{tikzpicture}[x=1cm,y=1cm,>=Stealth]
\draw[line width=3.4pt,black!22] (0.000,0) -- (0.043,0);
\draw[line width=3.4pt,black!22] (0.059,0) -- (0.071,0);
\draw[line width=3.4pt,black!22] (0.098,0) -- (0.118,0);
\draw[line width=3.4pt,black!22] (0.164,0) -- (0.197,0);
\draw[line width=3.4pt,black!22] (0.273,0) -- (0.328,0);
\draw[line width=3.4pt,black!22] (0.456,0) -- (0.547,0);
\draw[line width=3.4pt,black!22] (0.759,0) -- (0.912,0);
\draw[line width=3.4pt,black!22] (1.266,0) -- (1.520,0);
\draw[line width=3.4pt,black!22] (2.109,0) -- (2.533,0);
\draw[line width=3.4pt,black!22] (3.515,0) -- (4.221,0);
\draw[line width=3.4pt,black!22] (5.859,0) -- (7.035,0);
\draw[line width=3.4pt,black!22] (9.765,0) -- (11.112,0);
\draw[thin] plot coordinates {(7.035,0.000) (7.075,0.000) (7.114,0.000) (7.154,0.008) (7.193,0.035) (7.233,0.083) (7.272,0.147) (7.312,0.221) (7.352,0.299) (7.391,0.378) (7.431,0.456) (7.470,0.530) (7.510,0.601) (7.549,0.668) (7.589,0.730) (7.628,0.788) (7.668,0.842) (7.708,0.891) (7.747,0.937) (7.787,0.978) (7.826,1.017) (7.866,1.052) (7.905,1.083) (7.945,1.112) (7.985,1.138) (8.024,1.161) (8.064,1.181) (8.103,1.199) (8.143,1.214) (8.182,1.228) (8.222,1.238) (8.262,1.247) (8.301,1.253) (8.341,1.258) (8.380,1.260) (8.420,1.260) (8.459,1.258) (8.499,1.253) (8.538,1.247) (8.578,1.238) (8.618,1.228) (8.657,1.214) (8.697,1.199) (8.736,1.181) (8.776,1.161) (8.815,1.138) (8.855,1.112) (8.895,1.083) (8.934,1.052) (8.974,1.017) (9.013,0.978) (9.053,0.937) (9.092,0.891) (9.132,0.842) (9.172,0.788) (9.211,0.730) (9.251,0.668) (9.290,0.601) (9.330,0.530) (9.369,0.456) (9.409,0.378) (9.448,0.299) (9.488,0.221) (9.528,0.147) (9.567,0.083) (9.607,0.035) (9.646,0.008) (9.686,0.000) (9.725,0.000) (9.765,0.000)};
\draw[thin] (7.035,0) -- (9.765,0);
\draw[thin] plot coordinates {(4.221,0.000) (4.245,-0.000) (4.268,-0.000) (4.292,-0.006) (4.316,-0.026) (4.340,-0.061) (4.363,-0.108) (4.387,-0.162) (4.411,-0.220) (4.435,-0.278) (4.458,-0.335) (4.482,-0.390) (4.506,-0.442) (4.530,-0.491) (4.553,-0.537) (4.577,-0.580) (4.601,-0.619) (4.625,-0.656) (4.648,-0.689) (4.672,-0.720) (4.696,-0.748) (4.720,-0.774) (4.743,-0.797) (4.767,-0.818) (4.791,-0.837) (4.814,-0.854) (4.838,-0.869) (4.862,-0.883) (4.886,-0.894) (4.909,-0.904) (4.933,-0.911) (4.957,-0.918) (4.981,-0.923) (5.004,-0.926) (5.028,-0.927) (5.052,-0.927) (5.076,-0.926) (5.099,-0.923) (5.123,-0.918) (5.147,-0.911) (5.171,-0.904) (5.194,-0.894) (5.218,-0.883) (5.242,-0.869) (5.266,-0.854) (5.289,-0.837) (5.313,-0.818) (5.337,-0.797) (5.360,-0.774) (5.384,-0.748) (5.408,-0.720) (5.432,-0.689) (5.455,-0.656) (5.479,-0.619) (5.503,-0.580) (5.527,-0.537) (5.550,-0.491) (5.574,-0.442) (5.598,-0.390) (5.622,-0.335) (5.645,-0.278) (5.669,-0.220) (5.693,-0.162) (5.717,-0.108) (5.740,-0.061) (5.764,-0.026) (5.788,-0.006) (5.812,-0.000) (5.835,-0.000) (5.859,0.000)};
\draw[thin] (4.221,0) -- (5.859,0);
\draw[thin] plot coordinates {(2.533,0.000) (2.547,0.000) (2.561,0.000) (2.575,0.005) (2.590,0.019) (2.604,0.045) (2.618,0.080) (2.632,0.120) (2.647,0.162) (2.661,0.205) (2.675,0.247) (2.689,0.287) (2.704,0.326) (2.718,0.362) (2.732,0.395) (2.746,0.427) (2.760,0.456) (2.775,0.483) (2.789,0.507) (2.803,0.530) (2.817,0.551) (2.832,0.570) (2.846,0.587) (2.860,0.602) (2.874,0.616) (2.889,0.629) (2.903,0.640) (2.917,0.650) (2.931,0.658) (2.946,0.665) (2.960,0.671) (2.974,0.676) (2.988,0.679) (3.003,0.681) (3.017,0.682) (3.031,0.682) (3.045,0.681) (3.060,0.679) (3.074,0.676) (3.088,0.671) (3.102,0.665) (3.117,0.658) (3.131,0.650) (3.145,0.640) (3.159,0.629) (3.174,0.616) (3.188,0.602) (3.202,0.587) (3.216,0.570) (3.231,0.551) (3.245,0.530) (3.259,0.507) (3.273,0.483) (3.288,0.456) (3.302,0.427) (3.316,0.395) (3.330,0.362) (3.344,0.326) (3.359,0.287) (3.373,0.247) (3.387,0.205) (3.401,0.162) (3.416,0.120) (3.430,0.080) (3.444,0.045) (3.458,0.019) (3.473,0.005) (3.487,0.000) (3.501,0.000) (3.515,0.000)};
\draw[thin] (2.533,0) -- (3.515,0);
\draw[thin] plot coordinates {(1.520,0.000) (1.528,-0.000) (1.537,-0.000) (1.545,-0.003) (1.554,-0.014) (1.562,-0.033) (1.571,-0.059) (1.579,-0.088) (1.588,-0.119) (1.596,-0.151) (1.605,-0.182) (1.614,-0.211) (1.622,-0.240) (1.631,-0.266) (1.639,-0.291) (1.648,-0.314) (1.656,-0.336) (1.665,-0.355) (1.673,-0.373) (1.682,-0.390) (1.690,-0.405) (1.699,-0.419) (1.708,-0.432) (1.716,-0.443) (1.725,-0.454) (1.733,-0.463) (1.742,-0.471) (1.750,-0.478) (1.759,-0.484) (1.767,-0.489) (1.776,-0.494) (1.784,-0.497) (1.793,-0.500) (1.802,-0.501) (1.810,-0.502) (1.819,-0.502) (1.827,-0.501) (1.836,-0.500) (1.844,-0.497) (1.853,-0.494) (1.861,-0.489) (1.870,-0.484) (1.878,-0.478) (1.887,-0.471) (1.896,-0.463) (1.904,-0.454) (1.913,-0.443) (1.921,-0.432) (1.930,-0.419) (1.938,-0.405) (1.947,-0.390) (1.955,-0.373) (1.964,-0.355) (1.973,-0.336) (1.981,-0.314) (1.990,-0.291) (1.998,-0.266) (2.007,-0.240) (2.015,-0.211) (2.024,-0.182) (2.032,-0.151) (2.041,-0.119) (2.049,-0.088) (2.058,-0.059) (2.067,-0.033) (2.075,-0.014) (2.084,-0.003) (2.092,-0.000) (2.101,-0.000) (2.109,0.000)};
\draw[thin] (1.520,0) -- (2.109,0);
\draw[thin] plot coordinates {(0.912,0.000) (0.917,0.000) (0.922,0.000) (0.927,0.002) (0.932,0.010) (0.937,0.024) (0.943,0.043) (0.948,0.065) (0.953,0.088) (0.958,0.111) (0.963,0.134) (0.968,0.156) (0.973,0.176) (0.978,0.196) (0.984,0.214) (0.989,0.231) (0.994,0.247) (0.999,0.262) (1.004,0.275) (1.009,0.287) (1.014,0.298) (1.019,0.309) (1.025,0.318) (1.030,0.326) (1.035,0.334) (1.040,0.341) (1.045,0.347) (1.050,0.352) (1.055,0.356) (1.060,0.360) (1.066,0.363) (1.071,0.366) (1.076,0.368) (1.081,0.369) (1.086,0.370) (1.091,0.370) (1.096,0.369) (1.101,0.368) (1.107,0.366) (1.112,0.363) (1.117,0.360) (1.122,0.356) (1.127,0.352) (1.132,0.347) (1.137,0.341) (1.142,0.334) (1.148,0.326) (1.153,0.318) (1.158,0.309) (1.163,0.298) (1.168,0.287) (1.173,0.275) (1.178,0.262) (1.184,0.247) (1.189,0.231) (1.194,0.214) (1.199,0.196) (1.204,0.176) (1.209,0.156) (1.214,0.134) (1.219,0.111) (1.225,0.088) (1.230,0.065) (1.235,0.043) (1.240,0.024) (1.245,0.010) (1.250,0.002) (1.255,0.000) (1.260,0.000) (1.266,0.000)};
\draw[thin] (0.912,0) -- (1.266,0);
\draw[thin] plot coordinates {(0.547,-0.000) (0.550,-0.000) (0.553,-0.000) (0.556,-0.002) (0.559,-0.008) (0.562,-0.018) (0.566,-0.032) (0.569,-0.048) (0.572,-0.065) (0.575,-0.082) (0.578,-0.098) (0.581,-0.115) (0.584,-0.130) (0.587,-0.144) (0.590,-0.158) (0.593,-0.170) (0.596,-0.182) (0.599,-0.192) (0.602,-0.202) (0.605,-0.211) (0.609,-0.220) (0.612,-0.227) (0.615,-0.234) (0.618,-0.240) (0.621,-0.246) (0.624,-0.251) (0.627,-0.255) (0.630,-0.259) (0.633,-0.262) (0.636,-0.265) (0.639,-0.267) (0.642,-0.269) (0.645,-0.271) (0.649,-0.272) (0.652,-0.272) (0.655,-0.272) (0.658,-0.272) (0.661,-0.271) (0.664,-0.269) (0.667,-0.267) (0.670,-0.265) (0.673,-0.262) (0.676,-0.259) (0.679,-0.255) (0.682,-0.251) (0.685,-0.246) (0.689,-0.240) (0.692,-0.234) (0.695,-0.227) (0.698,-0.220) (0.701,-0.211) (0.704,-0.202) (0.707,-0.192) (0.710,-0.182) (0.713,-0.170) (0.716,-0.158) (0.719,-0.144) (0.722,-0.130) (0.725,-0.115) (0.729,-0.098) (0.732,-0.082) (0.735,-0.065) (0.738,-0.048) (0.741,-0.032) (0.744,-0.018) (0.747,-0.008) (0.750,-0.002) (0.753,-0.000) (0.756,-0.000) (0.759,-0.000)};
\draw[thin] (0.547,0) -- (0.759,0);
\draw[thin] plot coordinates {(0.328,0.000) (0.330,0.000) (0.332,0.000) (0.334,0.001) (0.336,0.006) (0.337,0.013) (0.339,0.023) (0.341,0.035) (0.343,0.048) (0.345,0.060) (0.347,0.072) (0.349,0.084) (0.350,0.096) (0.352,0.106) (0.354,0.116) (0.356,0.125) (0.358,0.134) (0.360,0.142) (0.361,0.149) (0.363,0.156) (0.365,0.162) (0.367,0.167) (0.369,0.172) (0.371,0.177) (0.373,0.181) (0.374,0.185) (0.376,0.188) (0.378,0.191) (0.380,0.193) (0.382,0.195) (0.384,0.197) (0.385,0.198) (0.387,0.199) (0.389,0.200) (0.391,0.200) (0.393,0.200) (0.395,0.200) (0.397,0.199) (0.398,0.198) (0.400,0.197) (0.402,0.195) (0.404,0.193) (0.406,0.191) (0.408,0.188) (0.409,0.185) (0.411,0.181) (0.413,0.177) (0.415,0.172) (0.417,0.167) (0.419,0.162) (0.421,0.156) (0.422,0.149) (0.424,0.142) (0.426,0.134) (0.428,0.125) (0.430,0.116) (0.432,0.106) (0.433,0.096) (0.435,0.084) (0.437,0.072) (0.439,0.060) (0.441,0.048) (0.443,0.035) (0.445,0.023) (0.446,0.013) (0.448,0.006) (0.450,0.001) (0.452,0.000) (0.454,0.000) (0.456,0.000)};
\draw[thin] (0.328,0) -- (0.456,0);
\draw[thin] plot coordinates {(0.197,-0.000) (0.198,-0.000) (0.199,-0.000) (0.200,-0.001) (0.201,-0.004) (0.202,-0.010) (0.204,-0.017) (0.205,-0.026) (0.206,-0.035) (0.207,-0.044) (0.208,-0.053) (0.209,-0.062) (0.210,-0.070) (0.211,-0.078) (0.212,-0.085) (0.214,-0.092) (0.215,-0.098) (0.216,-0.104) (0.217,-0.110) (0.218,-0.114) (0.219,-0.119) (0.220,-0.123) (0.221,-0.127) (0.222,-0.130) (0.224,-0.133) (0.225,-0.136) (0.226,-0.138) (0.227,-0.140) (0.228,-0.142) (0.229,-0.144) (0.230,-0.145) (0.231,-0.146) (0.232,-0.147) (0.233,-0.147) (0.235,-0.147) (0.236,-0.147) (0.237,-0.147) (0.238,-0.147) (0.239,-0.146) (0.240,-0.145) (0.241,-0.144) (0.242,-0.142) (0.243,-0.140) (0.245,-0.138) (0.246,-0.136) (0.247,-0.133) (0.248,-0.130) (0.249,-0.127) (0.250,-0.123) (0.251,-0.119) (0.252,-0.114) (0.253,-0.110) (0.255,-0.104) (0.256,-0.098) (0.257,-0.092) (0.258,-0.085) (0.259,-0.078) (0.260,-0.070) (0.261,-0.062) (0.262,-0.053) (0.263,-0.044) (0.264,-0.035) (0.266,-0.026) (0.267,-0.017) (0.268,-0.010) (0.269,-0.004) (0.270,-0.001) (0.271,-0.000) (0.272,-0.000) (0.273,-0.000)};
\draw[thin] (0.197,0) -- (0.273,0);
\draw[thin] plot coordinates {(0.118,0.000) (0.119,0.000) (0.119,0.000) (0.120,0.001) (0.121,0.003) (0.121,0.007) (0.122,0.013) (0.123,0.019) (0.123,0.026) (0.124,0.033) (0.125,0.039) (0.125,0.046) (0.126,0.052) (0.127,0.058) (0.127,0.063) (0.128,0.068) (0.129,0.072) (0.129,0.077) (0.130,0.081) (0.131,0.084) (0.131,0.088) (0.132,0.091) (0.133,0.093) (0.133,0.096) (0.134,0.098) (0.135,0.100) (0.135,0.102) (0.136,0.103) (0.137,0.105) (0.137,0.106) (0.138,0.107) (0.139,0.107) (0.139,0.108) (0.140,0.108) (0.141,0.108) (0.141,0.108) (0.142,0.108) (0.143,0.108) (0.143,0.107) (0.144,0.107) (0.145,0.106) (0.145,0.105) (0.146,0.103) (0.147,0.102) (0.147,0.100) (0.148,0.098) (0.149,0.096) (0.149,0.093) (0.150,0.091) (0.151,0.088) (0.151,0.084) (0.152,0.081) (0.153,0.077) (0.153,0.072) (0.154,0.068) (0.155,0.063) (0.155,0.058) (0.156,0.052) (0.157,0.046) (0.157,0.039) (0.158,0.033) (0.159,0.026) (0.159,0.019) (0.160,0.013) (0.161,0.007) (0.161,0.003) (0.162,0.001) (0.163,0.000) (0.163,0.000) (0.164,0.000)};
\draw[thin] (0.118,0) -- (0.164,0);
\draw[thin] plot coordinates {(0.071,-0.000) (0.071,-0.000) (0.072,-0.000) (0.072,-0.001) (0.072,-0.002) (0.073,-0.005) (0.073,-0.009) (0.074,-0.014) (0.074,-0.019) (0.074,-0.024) (0.075,-0.029) (0.075,-0.034) (0.076,-0.038) (0.076,-0.042) (0.076,-0.046) (0.077,-0.050) (0.077,-0.053) (0.078,-0.056) (0.078,-0.059) (0.078,-0.062) (0.079,-0.064) (0.079,-0.067) (0.080,-0.069) (0.080,-0.070) (0.080,-0.072) (0.081,-0.074) (0.081,-0.075) (0.082,-0.076) (0.082,-0.077) (0.082,-0.078) (0.083,-0.079) (0.083,-0.079) (0.084,-0.079) (0.084,-0.080) (0.084,-0.080) (0.085,-0.080) (0.085,-0.080) (0.086,-0.079) (0.086,-0.079) (0.086,-0.079) (0.087,-0.078) (0.087,-0.077) (0.088,-0.076) (0.088,-0.075) (0.088,-0.074) (0.089,-0.072) (0.089,-0.070) (0.090,-0.069) (0.090,-0.067) (0.090,-0.064) (0.091,-0.062) (0.091,-0.059) (0.092,-0.056) (0.092,-0.053) (0.092,-0.050) (0.093,-0.046) (0.093,-0.042) (0.094,-0.038) (0.094,-0.034) (0.094,-0.029) (0.095,-0.024) (0.095,-0.019) (0.096,-0.014) (0.096,-0.009) (0.096,-0.005) (0.097,-0.002) (0.097,-0.001) (0.098,-0.000) (0.098,-0.000) (0.098,-0.000)};
\draw[thin] (0.071,0) -- (0.098,0);
\draw[thin] plot coordinates {(0.043,0.000) (0.043,0.000) (0.043,0.000) (0.043,0.000) (0.043,0.002) (0.044,0.004) (0.044,0.007) (0.044,0.010) (0.044,0.014) (0.045,0.018) (0.045,0.021) (0.045,0.025) (0.045,0.028) (0.046,0.031) (0.046,0.034) (0.046,0.037) (0.046,0.039) (0.047,0.042) (0.047,0.044) (0.047,0.046) (0.047,0.047) (0.048,0.049) (0.048,0.051) (0.048,0.052) (0.048,0.053) (0.049,0.054) (0.049,0.055) (0.049,0.056) (0.049,0.057) (0.049,0.057) (0.050,0.058) (0.050,0.058) (0.050,0.058) (0.050,0.059) (0.051,0.059) (0.051,0.059) (0.051,0.059) (0.051,0.058) (0.052,0.058) (0.052,0.058) (0.052,0.057) (0.052,0.057) (0.053,0.056) (0.053,0.055) (0.053,0.054) (0.053,0.053) (0.054,0.052) (0.054,0.051) (0.054,0.049) (0.054,0.047) (0.054,0.046) (0.055,0.044) (0.055,0.042) (0.055,0.039) (0.055,0.037) (0.056,0.034) (0.056,0.031) (0.056,0.028) (0.056,0.025) (0.057,0.021) (0.057,0.018) (0.057,0.014) (0.057,0.010) (0.058,0.007) (0.058,0.004) (0.058,0.002) (0.058,0.000) (0.059,0.000) (0.059,0.000) (0.059,0.000)};
\draw[thin] (0.043,0) -- (0.059,0);
\draw[very thin,black!45] plot[smooth cycle] coordinates {(4.574,-0.278) (4.592,-0.208) (4.611,-0.184) (4.630,-0.169) (4.649,-0.158) (4.667,-0.149) (4.686,-0.142) (4.705,-0.137) (4.724,-0.133) (4.742,-0.129) (4.761,-0.126) (4.780,-0.123) (4.799,-0.121) (4.818,-0.119) (4.836,-0.118) (4.855,-0.116) (4.874,-0.115) (4.893,-0.115) (4.912,-0.114) (4.930,-0.114) (4.949,-0.114) (4.968,-0.114) (4.987,-0.114) (5.006,-0.114) (5.024,-0.115) (5.043,-0.115) (5.062,-0.116) (5.081,-0.117) (5.100,-0.119) (5.118,-0.120) (5.137,-0.122) (5.156,-0.124) (5.175,-0.126) (5.194,-0.128) (5.212,-0.131) (5.231,-0.134) (5.250,-0.137) (5.269,-0.141) (5.288,-0.146) (5.306,-0.151) (5.325,-0.157) (5.344,-0.164) (5.363,-0.173) (5.381,-0.184) (5.400,-0.197) (5.419,-0.215) (5.438,-0.242) (5.457,-0.322) (5.438,-0.439) (5.419,-0.491) (5.400,-0.533) (5.381,-0.568) (5.363,-0.599) (5.344,-0.626) (5.325,-0.651) (5.306,-0.673) (5.288,-0.693) (5.269,-0.711) (5.250,-0.727) (5.231,-0.742) (5.212,-0.755) (5.194,-0.766) (5.175,-0.776) (5.156,-0.785) (5.137,-0.793) (5.118,-0.799) (5.100,-0.804) (5.081,-0.808) (5.062,-0.810) (5.043,-0.812) (5.024,-0.812) (5.006,-0.812) (4.987,-0.810) (4.968,-0.807) (4.949,-0.802) (4.930,-0.797) (4.912,-0.790) (4.893,-0.782) (4.874,-0.773) (4.855,-0.763) (4.836,-0.751) (4.818,-0.737) (4.799,-0.723) (4.780,-0.706) (4.761,-0.688) (4.742,-0.668) (4.724,-0.646) (4.705,-0.622) (4.686,-0.595) (4.667,-0.565) (4.649,-0.532) (4.630,-0.495) (4.611,-0.452) (4.592,-0.398)};
\draw[very thin,black!45] plot[smooth cycle] coordinates {(4.726,-0.385) (4.738,-0.332) (4.749,-0.313) (4.761,-0.299) (4.773,-0.289) (4.784,-0.280) (4.796,-0.273) (4.808,-0.267) (4.819,-0.263) (4.831,-0.258) (4.843,-0.255) (4.855,-0.252) (4.866,-0.249) (4.878,-0.247) (4.890,-0.245) (4.902,-0.243) (4.913,-0.242) (4.925,-0.241) (4.937,-0.240) (4.948,-0.240) (4.960,-0.240) (4.972,-0.239) (4.984,-0.240) (4.995,-0.240) (5.007,-0.241) (5.019,-0.241) (5.031,-0.242) (5.042,-0.244) (5.054,-0.245) (5.066,-0.247) (5.077,-0.249) (5.089,-0.251) (5.101,-0.254) (5.112,-0.256) (5.124,-0.260) (5.136,-0.263) (5.148,-0.267) (5.159,-0.271) (5.171,-0.276) (5.183,-0.282) (5.195,-0.288) (5.206,-0.295) (5.218,-0.304) (5.230,-0.314) (5.241,-0.326) (5.253,-0.341) (5.265,-0.362) (5.277,-0.415) (5.265,-0.492) (5.253,-0.521) (5.241,-0.544) (5.230,-0.562) (5.218,-0.579) (5.206,-0.593) (5.195,-0.606) (5.183,-0.617) (5.171,-0.627) (5.159,-0.636) (5.148,-0.644) (5.136,-0.651) (5.124,-0.658) (5.112,-0.664) (5.101,-0.669) (5.089,-0.673) (5.077,-0.677) (5.066,-0.680) (5.054,-0.682) (5.042,-0.684) (5.031,-0.685) (5.019,-0.685) (5.007,-0.685) (4.995,-0.685) (4.984,-0.683) (4.972,-0.681) (4.960,-0.679) (4.948,-0.676) (4.937,-0.672) (4.925,-0.668) (4.913,-0.663) (4.902,-0.657) (4.890,-0.651) (4.878,-0.644) (4.866,-0.636) (4.855,-0.627) (4.843,-0.617) (4.831,-0.607) (4.819,-0.595) (4.808,-0.582) (4.796,-0.568) (4.784,-0.553) (4.773,-0.535) (4.761,-0.514) (4.749,-0.490) (4.738,-0.460)};
\draw[very thin,black!45] plot[smooth cycle] coordinates {(2.752,0.224) (2.764,0.302) (2.775,0.341) (2.786,0.372) (2.798,0.399) (2.809,0.423) (2.820,0.445) (2.832,0.464) (2.843,0.481) (2.854,0.497) (2.866,0.511) (2.877,0.524) (2.888,0.536) (2.899,0.547) (2.911,0.556) (2.922,0.564) (2.933,0.572) (2.945,0.578) (2.956,0.584) (2.967,0.588) (2.979,0.592) (2.990,0.595) (3.001,0.596) (3.013,0.597) (3.024,0.598) (3.035,0.597) (3.047,0.596) (3.058,0.593) (3.069,0.590) (3.081,0.586) (3.092,0.581) (3.103,0.575) (3.115,0.568) (3.126,0.560) (3.137,0.551) (3.149,0.541) (3.160,0.530) (3.171,0.518) (3.182,0.504) (3.194,0.489) (3.205,0.472) (3.217,0.453) (3.228,0.433) (3.239,0.409) (3.250,0.383) (3.262,0.352) (3.273,0.313) (3.284,0.233) (3.273,0.170) (3.262,0.151) (3.250,0.138) (3.239,0.129) (3.228,0.122) (3.217,0.116) (3.205,0.111) (3.194,0.107) (3.182,0.104) (3.171,0.101) (3.160,0.098) (3.149,0.096) (3.137,0.094) (3.126,0.092) (3.115,0.091) (3.103,0.090) (3.092,0.088) (3.081,0.088) (3.069,0.087) (3.058,0.086) (3.047,0.086) (3.035,0.085) (3.024,0.085) (3.013,0.085) (3.001,0.085) (2.990,0.085) (2.979,0.085) (2.967,0.085) (2.956,0.086) (2.945,0.086) (2.933,0.087) (2.922,0.088) (2.911,0.089) (2.899,0.091) (2.888,0.092) (2.877,0.094) (2.866,0.096) (2.854,0.099) (2.843,0.102) (2.832,0.106) (2.820,0.110) (2.809,0.115) (2.798,0.122) (2.786,0.130) (2.775,0.142) (2.764,0.160)};
\draw[very thin,black!45] plot[smooth cycle] coordinates {(2.850,0.298) (2.857,0.346) (2.864,0.368) (2.871,0.385) (2.878,0.400) (2.885,0.413) (2.892,0.424) (2.899,0.434) (2.906,0.443) (2.913,0.451) (2.920,0.458) (2.927,0.465) (2.934,0.471) (2.941,0.476) (2.948,0.481) (2.955,0.486) (2.962,0.490) (2.969,0.493) (2.976,0.496) (2.983,0.498) (2.990,0.500) (2.997,0.502) (3.005,0.503) (3.012,0.503) (3.019,0.504) (3.026,0.503) (3.033,0.503) (3.040,0.501) (3.047,0.500) (3.054,0.498) (3.061,0.495) (3.068,0.492) (3.075,0.489) (3.082,0.485) (3.089,0.481) (3.096,0.475) (3.103,0.470) (3.110,0.464) (3.117,0.457) (3.124,0.449) (3.131,0.440) (3.138,0.431) (3.145,0.420) (3.152,0.408) (3.159,0.394) (3.166,0.377) (3.173,0.356) (3.181,0.310) (3.173,0.261) (3.166,0.246) (3.159,0.235) (3.152,0.227) (3.145,0.220) (3.138,0.214) (3.131,0.209) (3.124,0.205) (3.117,0.201) (3.110,0.198) (3.103,0.195) (3.096,0.192) (3.089,0.190) (3.082,0.188) (3.075,0.186) (3.068,0.185) (3.061,0.183) (3.054,0.182) (3.047,0.181) (3.040,0.180) (3.033,0.180) (3.026,0.179) (3.019,0.179) (3.012,0.179) (3.005,0.179) (2.997,0.179) (2.990,0.179) (2.983,0.180) (2.976,0.180) (2.969,0.181) (2.962,0.182) (2.955,0.183) (2.948,0.185) (2.941,0.186) (2.934,0.188) (2.927,0.191) (2.920,0.193) (2.913,0.196) (2.906,0.199) (2.899,0.203) (2.892,0.208) (2.885,0.213) (2.878,0.220) (2.871,0.228) (2.864,0.238) (2.857,0.253)};
\draw[line width=1.15pt,blue!75!black] plot coordinates {(0.000,0.000) (0.043,0.000) (0.043,0.000) (0.043,0.000) (0.043,0.000) (0.043,0.000) (0.043,0.002) (0.044,0.004) (0.044,0.007) (0.044,0.010) (0.044,0.014) (0.045,0.018) (0.045,0.021) (0.045,0.025) (0.045,0.028) (0.046,0.031) (0.046,0.034) (0.046,0.037) (0.046,0.039) (0.047,0.042) (0.047,0.044) (0.047,0.046) (0.047,0.047) (0.048,0.049) (0.048,0.051) (0.048,0.052) (0.048,0.053) (0.049,0.054) (0.049,0.055) (0.049,0.056) (0.049,0.057) (0.049,0.057) (0.050,0.058) (0.050,0.058) (0.050,0.058) (0.050,0.059) (0.051,0.059) (0.051,0.059) (0.051,0.059) (0.051,0.058) (0.052,0.058) (0.052,0.058) (0.052,0.057) (0.052,0.057) (0.053,0.056) (0.053,0.055) (0.053,0.054) (0.053,0.053) (0.054,0.052) (0.054,0.051) (0.054,0.049) (0.054,0.047) (0.054,0.046) (0.055,0.044) (0.055,0.042) (0.055,0.039) (0.055,0.037) (0.056,0.034) (0.056,0.031) (0.056,0.028) (0.056,0.025) (0.057,0.021) (0.057,0.018) (0.057,0.014) (0.057,0.010) (0.058,0.007) (0.058,0.004) (0.058,0.002) (0.058,0.000) (0.059,0.000) (0.059,0.000) (0.059,0.000) (0.071,0.000) (0.098,0.000) (0.118,0.000) (0.118,0.000) (0.119,0.000) (0.119,0.000) (0.120,0.001) (0.121,0.003) (0.121,0.007) (0.122,0.013) (0.123,0.019) (0.123,0.026) (0.124,0.033) (0.125,0.039) (0.125,0.046) (0.126,0.052) (0.127,0.058) (0.127,0.063) (0.128,0.068) (0.129,0.072) (0.129,0.077) (0.130,0.081) (0.131,0.084) (0.131,0.088) (0.132,0.091) (0.133,0.093) (0.133,0.096) (0.134,0.098) (0.135,0.100) (0.135,0.102) (0.136,0.103) (0.137,0.105) (0.137,0.106) (0.138,0.107) (0.139,0.107) (0.139,0.108) (0.140,0.108) (0.141,0.108) (0.141,0.108) (0.142,0.108) (0.143,0.108) (0.143,0.107) (0.144,0.107) (0.145,0.106) (0.145,0.105) (0.146,0.103) (0.147,0.102) (0.147,0.100) (0.148,0.098) (0.149,0.096) (0.149,0.093) (0.150,0.091) (0.151,0.088) (0.151,0.084) (0.152,0.081) (0.153,0.077) (0.153,0.072) (0.154,0.068) (0.155,0.063) (0.155,0.058) (0.156,0.052) (0.157,0.046) (0.157,0.039) (0.158,0.033) (0.159,0.026) (0.159,0.019) (0.160,0.013) (0.161,0.007) (0.161,0.003) (0.162,0.001) (0.163,0.000) (0.163,0.000) (0.164,0.000) (0.197,0.000) (0.273,0.000) (0.328,0.000) (0.328,0.000) (0.330,0.000) (0.332,0.000) (0.334,0.001) (0.336,0.006) (0.337,0.013) (0.339,0.023) (0.341,0.035) (0.343,0.048) (0.345,0.060) (0.347,0.072) (0.349,0.084) (0.350,0.096) (0.352,0.106) (0.354,0.116) (0.356,0.125) (0.358,0.134) (0.360,0.142) (0.361,0.149) (0.363,0.156) (0.365,0.162) (0.367,0.167) (0.369,0.172) (0.371,0.177) (0.373,0.181) (0.374,0.185) (0.376,0.188) (0.378,0.191) (0.380,0.193) (0.382,0.195) (0.384,0.197) (0.385,0.198) (0.387,0.199) (0.389,0.200) (0.391,0.200) (0.393,0.200) (0.395,0.200) (0.397,0.199) (0.398,0.198) (0.400,0.197) (0.402,0.195) (0.404,0.193) (0.406,0.191) (0.408,0.188) (0.409,0.185) (0.411,0.181) (0.413,0.177) (0.415,0.172) (0.417,0.167) (0.419,0.162) (0.421,0.156) (0.422,0.149) (0.424,0.142) (0.426,0.134) (0.428,0.125) (0.430,0.116) (0.432,0.106) (0.433,0.096) (0.435,0.084) (0.437,0.072) (0.439,0.060) (0.441,0.048) (0.443,0.035) (0.445,0.023) (0.446,0.013) (0.448,0.006) (0.450,0.001) (0.452,0.000) (0.454,0.000) (0.456,0.000) (0.547,0.000) (0.547,-0.000) (0.550,-0.000) (0.553,-0.000) (0.556,-0.002) (0.559,-0.008) (0.562,-0.018) (0.566,-0.032) (0.569,-0.048) (0.572,-0.065) (0.575,-0.082) (0.578,-0.098) (0.581,-0.115) (0.584,-0.130) (0.587,-0.144) (0.590,-0.158) (0.593,-0.170) (0.596,-0.182) (0.599,-0.192) (0.602,-0.202) (0.605,-0.211) (0.609,-0.220) (0.612,-0.227) (0.615,-0.234) (0.618,-0.240) (0.621,-0.246) (0.624,-0.251) (0.627,-0.255) (0.630,-0.259) (0.633,-0.262) (0.636,-0.265) (0.639,-0.267) (0.642,-0.269) (0.645,-0.271) (0.649,-0.272) (0.652,-0.272) (0.655,-0.272) (0.658,-0.272) (0.661,-0.271) (0.664,-0.269) (0.667,-0.267) (0.670,-0.265) (0.673,-0.262) (0.676,-0.259) (0.679,-0.255) (0.682,-0.251) (0.685,-0.246) (0.689,-0.240) (0.692,-0.234) (0.695,-0.227) (0.698,-0.220) (0.701,-0.211) (0.704,-0.202) (0.707,-0.192) (0.710,-0.182) (0.713,-0.170) (0.716,-0.158) (0.719,-0.144) (0.722,-0.130) (0.725,-0.115) (0.729,-0.098) (0.732,-0.082) (0.735,-0.065) (0.738,-0.048) (0.741,-0.032) (0.744,-0.018) (0.747,-0.008) (0.750,-0.002) (0.753,-0.000) (0.756,-0.000) (0.759,-0.000) (0.912,0.000) (1.266,0.000) (1.520,0.000) (1.520,0.000) (1.528,-0.000) (1.537,-0.000) (1.545,-0.003) (1.554,-0.014) (1.562,-0.033) (1.571,-0.059) (1.579,-0.088) (1.588,-0.119) (1.596,-0.151) (1.605,-0.182) (1.614,-0.211) (1.622,-0.240) (1.631,-0.266) (1.639,-0.291) (1.648,-0.314) (1.656,-0.336) (1.665,-0.355) (1.673,-0.373) (1.682,-0.390) (1.690,-0.405) (1.699,-0.419) (1.708,-0.432) (1.716,-0.443) (1.725,-0.454) (1.733,-0.463) (1.742,-0.471) (1.750,-0.478) (1.759,-0.484) (1.767,-0.489) (1.776,-0.494) (1.784,-0.497) (1.793,-0.500) (1.802,-0.501) (1.810,-0.502) (1.819,-0.502) (1.827,-0.501) (1.836,-0.500) (1.844,-0.497) (1.853,-0.494) (1.861,-0.489) (1.870,-0.484) (1.878,-0.478) (1.887,-0.471) (1.896,-0.463) (1.904,-0.454) (1.913,-0.443) (1.921,-0.432) (1.930,-0.419) (1.938,-0.405) (1.947,-0.390) (1.955,-0.373) (1.964,-0.355) (1.973,-0.336) (1.981,-0.314) (1.990,-0.291) (1.998,-0.266) (2.007,-0.240) (2.015,-0.211) (2.024,-0.182) (2.032,-0.151) (2.041,-0.119) (2.049,-0.088) (2.058,-0.059) (2.067,-0.033) (2.075,-0.014) (2.084,-0.003) (2.092,-0.000) (2.101,-0.000) (2.109,0.000) (2.533,0.000) (3.515,0.000) (4.221,0.000) (4.221,0.000) (4.245,-0.000) (4.268,-0.000) (4.292,-0.006) (4.316,-0.026) (4.340,-0.061) (4.363,-0.108) (4.387,-0.162) (4.411,-0.220) (4.435,-0.278) (4.458,-0.335) (4.482,-0.390) (4.506,-0.442) (4.530,-0.491) (4.553,-0.537) (4.577,-0.580) (4.601,-0.619) (4.625,-0.656) (4.648,-0.689) (4.672,-0.720) (4.696,-0.748) (4.720,-0.774) (4.743,-0.797) (4.767,-0.818) (4.791,-0.837) (4.814,-0.854) (4.838,-0.869) (4.862,-0.883) (4.886,-0.894) (4.909,-0.904) (4.933,-0.911) (4.957,-0.918) (4.981,-0.923) (5.004,-0.926) (5.028,-0.927) (5.052,-0.927) (5.076,-0.926) (5.099,-0.923) (5.123,-0.918) (5.147,-0.911) (5.171,-0.904) (5.194,-0.894) (5.218,-0.883) (5.242,-0.869) (5.266,-0.854) (5.289,-0.837) (5.313,-0.818) (5.337,-0.797) (5.360,-0.774) (5.384,-0.748) (5.408,-0.720) (5.432,-0.689) (5.455,-0.656) (5.479,-0.619) (5.503,-0.580) (5.527,-0.537) (5.550,-0.491) (5.574,-0.442) (5.598,-0.390) (5.622,-0.335) (5.645,-0.278) (5.669,-0.220) (5.693,-0.162) (5.717,-0.108) (5.740,-0.061) (5.764,-0.026) (5.788,-0.006) (5.812,-0.000) (5.835,-0.000) (5.859,0.000) (7.035,0.000) (8.400,0.000) (8.400,1.860)};
\fill (0.328,0) circle (0.8pt); \fill (0.456,0) circle (0.8pt);
\fill (0.547,0) circle (0.8pt); \fill (0.759,0) circle (0.8pt);
\fill (0.912,0) circle (0.8pt); \fill (1.266,0) circle (0.8pt);
\fill (1.520,0) circle (0.8pt); \fill (2.109,0) circle (0.8pt);
\fill (2.533,0) circle (0.8pt); \fill (3.515,0) circle (0.8pt);
\fill (4.221,0) circle (0.8pt); \fill (5.859,0) circle (0.8pt);
\fill (7.035,0) circle (0.8pt); \fill (9.765,0) circle (0.8pt);
\draw[dashed,thin] (8.400,-0.55) -- (8.400,2.25);
\draw[line width=1.15pt,blue!75!black] (8.400,0) -- (8.400,1.860);
\fill[blue!75!black] (8.400,1.860) circle (1.3pt);
\draw[fill=white] (8.400,0) circle (1.7pt); \draw[fill=white] (8.400,1.260) circle (1.7pt);
\node[font=\scriptsize,anchor=north west] at (8.430,-0.03) {$x_*$};
\node[font=\scriptsize,anchor=south] at (8.400,2.25) {$x=x_*$};
\fill (0,0) circle (1.1pt);
\node[font=\scriptsize,anchor=east] at (-0.05,0) {$0$};
\node[font=\scriptsize,anchor=west] at (11.112,0) {$y=0$};
\node[font=\scriptsize,anchor=south west] at (9.246,0.675) {$y=a(x)$};
\end{tikzpicture}
\caption{The zero level of $g_\nu$ near an accumulation point of the bumps (schematic). Over
each lobe of $a$ the axis and the graph are two orbits of $\mathcal X_\nu$ bounding an eye;
the thin closed curves are other level curves of $g_\nu$. The thick grey segments of the axis
consist of equilibria. A singular curve (blue) chooses one boundary orbit in each eye and
turns onto the invariant line $x=x_*$ at the saddle $(x_*,0)$.}
\label{fig:eyes}
\end{figure}

Our construction uses only the multipliers $\mu=(\nu,0)$. For them
$g_\nu:=g_{(\nu,0)}=\ell_\nu(x)\,y\,\bigl(y-a(x)\bigr)$, with $\ell_\nu=\nu\cdot w$, and
\[
\mathcal X_\nu=
\begin{cases}
-\ell_\nu a\,\partial_x & \text{on }\{y=0\},\\
\ell_\nu a\,(\partial_x+a'\partial_y) & \text{on }\{y=a(x)\},\\
-\ell_\nu'\,p\,\partial_y & \text{on }\{x=x_\nu\}.
\end{cases}
\]
Away from the line $x=x_\nu$, the zero level $\{y=0\}\cup\{y=a(x)\}$ is organised as follows
(Figure~\ref{fig:eyes}).
\begin{enumerate}
\item Over each gap of $a$, that is, an interval on which $a\equiv0$, the axis is a segment
of equilibria. There $g_\nu=\ell_\nu y^2$, the constraint is vacuous, and an abnormal
extremal may move along the segment in either direction while staying in $(\D^3)^\perp$.
\item Over each lobe of $a$, that is, each component of $\{a\ne0\}$, the axis arc and the
graph arc are two orbits of $\mathcal X_\nu$ joining the equilibria at the ends of the lobe,
traversed in opposite directions. Together they bound an \emph{eye}, like the region between
two consecutive saddles of the pendulum.
\item The eyes accumulate at the flat points $0,1,2$ of $a$.
\end{enumerate}
The line $x=x_\nu$ is invariant and meets the rest of the zero level at the equilibria
$(x_\nu,0)$ and $(x_\nu,a(x_\nu))$, which are hyperbolic saddles with eigenvalues
$\pm|a(x_\nu)\ell_\nu'|$ when $a(x_\nu)\ne0$.

On the axis $\det\operatorname{Hess}g_\nu=-(a'\ell_\nu+a\ell_\nu')^2\le0$, so no elliptic
points occur there. At the ends of the eyes $a=a'=0$, so these equilibria are nilpotent.
Since $a$ is flat there, the two orbits of an eye have infinite order of contact at its
ends, and no jet of $\D$ at such a point distinguishes it from a point of a gap. In the
present construction every branch point on the axis is of this nilpotent type: $a$ is a sum
of constant multiples of translated and rescaled copies of the nonnegative function chosen in the proof of
Lemma~\ref{lem:profiles}, with pairwise disjoint supports, so $a'=0$ wherever $a=0$.

\subsection*{The mechanism}
An abnormal extremal with multiplier $(\nu,0)$ is a concatenation of these orbits and
segments of equilibria, changing orbit only at equilibria. Crossing a string of eyes from
left to right, it chooses at each eye one of its two boundary orbits. This is Kupka's binary
cube \cite{Kupka}, realised by the abnormal dynamics. Following the graph arc of an eye
instead of its axis arc changes the vertical displacement by the strip vector of the
corresponding bump, up to sign (Lemma~\ref{lem:flux}). So the vertical displacements of
these extremals are subsums of a series of strip vectors, and
Sections~\ref{sec:moments}--\ref{sec:assemble} choose the eyes so that these displacements,
combined over finitely many loops, cover $\R^7$. The final turn onto the invariant line
$x=x_*$ at the equilibrium $(x_*,0)$ then reaches the prescribed base point.

\subsection*{What is essential}
Neither the infinite order of contact nor the nilpotent type of the branch points is
essential. In the transversal variant of the construction mentioned in the abstract, which
will be given in the longer version of the paper in preparation, every eye ends at
hyperbolic saddles, and near the axis the zero level is a chain of saddle connections, as for
the pendulum. There an extremal meets $(\D^3)^\perp$ only at isolated points: it arrives at a saddle
along one separatrix and may leave it along any other, and this is where the binary
choice is made. What cannot be dispensed with is the accumulation of infinitely many equilibria
along a chain of orbits of finite length. It turns the countable family of concatenations
available in a compact region into a Cantor set. For real-analytic $g_\mu$ this is impossible,
because the zero level is then locally a finite union of analytic arcs; this is the dynamical
reason why the mechanism needs a non-analytic coefficient. Equiregularity does not exclude it:
$(\D^2)^\perp$ and $(\D^3)^\perp$ are smooth subbundles of ranks $6$ and $4$, but the set of
branch points of the zero levels, which is what the abnormal dynamics sees, is not locally
finite.

\section*{AI-assisted tools disclosure}

In accordance with the Leiden Declaration on Artificial Intelligence and Mathematics
\cite{Alper}, we describe the automated tools used in this work.

\emph{Origin of the project.} The strategy of lifting the Lerario--Rizzi--Tiberio
example to higher dimensions, so that its singular curves lift to singular curves while
the lifted distribution becomes equiregular, was proposed by the second author. To
understand the geometric origin of the resulting example better, he also posed the
questions of whether its Tanaka symbol is constant and how its abnormal extremals
behave in the cotangent bundle, which are answered in Appendices~\ref{app:symbol}
and~\ref{app:cotangent}. He further suggested pushing the coefficient $a$ into
structure functions of positive weight, which led to the  counterexample in dimension twelve
with constant Tanaka symbol announced in Appendix~\ref{app:symbol}; details will
appear in \cite{ZW12}.

\emph{Large language models.} The mathematical implementation of this strategy was
generated in September and October 2026 in interactive sessions with Anthropic's
Claude [Claude Opus 5 and Claude Opus 5.5], used through the Claude desktop application. This includes the
explicit structure \eqref{eq:frame-intro}, the first versions of the proofs of
Theorems~\ref{thm:A}--\ref{thm:C}, in particular the formulation and solution of the
moment problem of Sections~\ref{sec:moments}--\ref{sec:assemble}, the statements in
Appendices~\ref{app:symbol} and~\ref{app:cotangent}, the explicit dimension-twelve
structure announced in Appendix~\ref{app:symbol}, and first drafts of most of the text, together
with revisions made in response to reviews. OpenAI's ChatGPT [model 6~Astra] was used
at a later stage to review drafts against the reports of the
formalisation and propose revisions.

\emph{Human contribution and verification.} Both authors checked all arguments line
by line. The revision of the AI-generated draft was led by the first author, and the
verification was supervised by the second author. The revision concentrated on the
technically most demanding, analytic part of the argument:
Sections~\ref{sec:moments}--\ref{sec:assemble} replace a single section of the
AI-generated draft, with some arguments simplified and others expanded. The
introduction, Sections~\ref{sec:lifts} and~\ref{sec:structure} and the appendices were
checked with the same care but needed fewer changes, because the AI-generated draft
already presented them in a simpler and more easily verifiable form.

\emph{Formal verification.} The first author formalised
Theorems~\ref{thm:A}--\ref{thm:C} in the proof assistant Lean~4 (version 4.30.0, with
Mathlib v4.30.0). The Lean formalisation was developed with AI assistance from OpenAI’s Codex, under the first author’s supervision. Codex assisted with drafting and revising Lean definitions, theorem statements and proof scripts, and debugging them using local compiler feedback.
The formalisation
uses the $L^2$ control formulation of Lemma~\ref{lem:endpoint} on $M=\{|y|<1\}$ and
covers the complete proofs of the three theorems, including both directions of
Lemma~\ref{lem:sing}, the graph-strip identity, the solution of the moment problem and
the construction of the coefficient; the identification $M\cong\R^9$ is proved
separately. Appendices~\ref{app:symbol} and~\ref{app:cotangent} are not formalised. The
development contains no unproved statements and uses only the standard axioms
\texttt{propext}, \texttt{Classical.choice} and \texttt{Quot.sound}. Some auxiliary
choices, such as the profile functions and several constants, differ from those in the
text. A successful build certifies the formal statements; that they correspond to the
statements of this paper has been checked by the authors. The formalisation
detected two defects in an earlier draft, an unrestricted uniformity claim in the
estimates for pairs of bumps and an unrestricted lifting clause, both of which have been
corrected. The sources, with build instructions, are available at \url{https://doi.org/10.5281/zenodo.23310061}.

\emph{Responsibility.} The authors take full and exclusive responsibility for the
correctness of the results, for the contents of the paper and for its citations.

\end{document}